\documentclass[11pt]{amsart} 
\usepackage[dvips]{graphicx}
\usepackage{amssymb,amsfonts,amsmath}
\newtheorem{theo}{Theorem}[section]
\newtheorem{prop}[theo]{Proposition}
\newtheorem{lemma}[theo]{Lemma}
\newtheorem{coro}[theo]{Corollary}

\newcommand{\Part}{{\mathcal P}}
\newcommand{\cross}{\mathrm {Cross}}
\newcommand{\K}{\mathcal{K}}
\newcommand{\sm}{\setminus}
\newcommand{\W}{\mathcal{W}}
\newcommand{\AP}{\mathcal{AP}}
\newcommand{\Sa}{\mathcal{S}}
\newcommand{\eps}{\varepsilon} 
\newcommand{\Hy}{\mathcal{H}}
\newcommand{\Pa}{{\mathcal P}}
\newcommand{\G}{{\mathcal G}}
\newcommand{\F}{{\mathcal F}}
\newcommand{\T}{{\mathcal T}}
\newcommand{\V}{\mathcal{V}}
\newcommand{\R}{\mathcal{R}}
\newcommand{\Q}{{\mathcal Q}}
\def\proof{\noindent{\bf Proof.}\ }

\def\COMMENT#1{}
\def\TASK#1{}
\def\noproof{{\unskip\nobreak\hfill\penalty50\hskip2em\hbox{}\nobreak\hfill%
       $\square$\parfillskip=0pt\finalhyphendemerits=0\par}\goodbreak}
\def\endproof{\noproof\bigskip}
\newdimen\margin   
\def\textno#1&#2\par{%
   \margin=\hsize
   \advance\margin by -4\parindent
          \setbox1=\hbox{\sl#1}%
   \ifdim\wd1 < \margin
      $$\box1\eqno#2$$%
   \else
      \bigbreak
      \hbox to \hsize{\indent$\vcenter{\advance\hsize by -3\parindent
      \sl\noindent#1}\hfil#2$}%
      \bigbreak
   \fi}
\def\proof{\removelastskip\penalty55\medskip\noindent{\bf Proof. }}
\def\enddiscard{}
\long\def\discard#1\enddiscard{}

\title{Hamilton $\ell$-cycles in uniform hypergraphs}
\author{Daniela K\"uhn \and Richard Mycroft \and Deryk Osthus}
\date{}
\thanks {D.~K\"uhn was partially supported by the EPSRC, grant no.~EP/D50564X/1.
D.~Osthus was partially supported by the EPSRC, grant no.~EP/E02162X/1.}

\begin{document}

\vspace*{-0.8cm}
\begin{abstract}
We say that a $k$-uniform hypergraph~$C$ is an $\ell$-\textit{cycle} if there 
exists a cyclic ordering of the vertices of~$C$ such that every edge of~$C$ 
consists of~$k$ consecutive vertices and such that every pair of {consecutive} edges 
(in the natural ordering of the edges) intersects in precisely $\ell$ vertices. 
We prove that if $1 \leq \ell < k$ and $k-\ell$ does not divide $k$ then any 
$k$-uniform hypergraph on $n$ vertices with minimum 
degree at least $\frac{n}{\lceil \frac{k}{k-\ell}\rceil (k-\ell)}+o(n)$ contains 
a Hamilton $\ell$-cycle. This confirms a conjecture of H\`an and Schacht. 
Together with results of R\"odl, Ruci\'nski and Szemer\'edi, our result 
asymptotically determines the minimum degree which forces an $\ell$-cycle for 
any $\ell$ with $1 \leq \ell <k$.  
\end{abstract}

\maketitle

\vspace*{-0.6cm}
\section{Introduction}\label{intro}

A $k$-\textit{graph}~$\Hy$ (also known as a 
$k$-uniform hypergraph), consists of a set of \textit{vertices}~$V(\Hy)$ and a 
set of \textit{edges}~$E(\Hy) \subseteq \{X \subseteq V(\Hy): |X| = k\}$, so 
that each edge of $\Hy$ consists of $k$ vertices. Let $\Hy$ be a 
$k$-graph, and let $A$ be a set of $k-1$ vertices of $\Hy$. Then 
the \textit{degree} of~$A$, denoted $d_\Hy(A)$, is the number of edges of $\Hy$ 
which contain $A$ as a subset. The \textit{minimum} \textit{degree} 
$\delta(\Hy)$ of $\Hy$ is then the minimum value of $d_\Hy(A)$ taken over all 
sets $A$ of $k-1$ vertices of $\Hy$.

We say that a $k$-graph~$C$ is an $\ell$-\textit{cycle} if there exists a cyclic 
ordering of the vertices of~$C$ such that every edge of~$C$ consists of~$k$ 
consecutive vertices and such that every pair of {consecutive} edges (in the natural 
ordering of the edges) intersects in precisely $\ell$ vertices. We say that a 
$k$-graph~$\Hy$ contains a \textit{Hamilton} $\ell$-\textit{cycle} if it 
contains a spanning sub-$k$-graph which is an $\ell$-cycle. Note that if a 
$k$-graph~$\Hy$ on~$n$ vertices contains a Hamilton $\ell$-cycle then $(k-\ell) 
| n$, since every edge {of the cycle} contains exactly $k-\ell$ vertices which were not 
contained in the previous edge.

We shall give an asymptotic solution to the question of what minimum degree will 
guarantee that a $k$-graph~$\Hy$ on~$n$ vertices contains a Hamilton 
$\ell$-cycle. This can be viewed as a generalisation of Dirac's theorem 
\cite{DIRAC}, which states that any graph (i.e. 2-graph) with $n \geq 3$ 
vertices and of minimum degree 
at least~$n/2$ contains a Hamilton cycle. 

In~\cite{RRS} and~\cite{RRS2}, R\"odl, Ruci\'nski and Szemer\'edi proved the 
following theorem for $\ell = k-1$; the other cases follow, since if 
$(k-\ell)|n$ 
then any $(k-1)$-cycle of order $n$ contains an $\ell$-cycle on the same 
vertices.

\begin{theo} \label{rrs}
For all~$k\geq 3$, $1 \leq \ell \leq k-1$ and any~$\eta > 0$ there exists~$n_0$ 
so that if~$n > n_0$ and~$(k-\ell)|n$ then any $k$-graph~$\Hy$ on~$n$ vertices 
with $\delta(\Hy) \geq \left(\frac{1}{2}+\eta\right) n$ contains a Hamilton 
$\ell$-cycle.
\end{theo}

This proved a conjecture of Katona and Kierstead~\cite{KK}.
Proposition~\ref{rrsbestpos} shows that Theorem~\ref{rrs} is best possible up to 
the error term~$\eta n$ if~$(k-\ell) | k$. This then raises the natural question of what 
minimum degree guarantees a Hamilton $\ell$-cycle if $(k-\ell) \nmid k$. 
In~\cite{KO}, K\"uhn and Osthus showed that any 3-graph $\Hy$ on $n$ vertices 
with $n$ even and $\delta(\Hy)\geq(\frac{1}{4}+o(1))n$ contains a Hamilton 
1-cycle. Keevash, K\"uhn, Mycroft and Osthus ~\cite{KKMO} extended this 
result to $k$-graphs, showing that any $k$-graph $\Hy$ on $n$ 
vertices with $(k-1)|n$ and  $\delta(\Hy)\geq(\frac{1}{2k-2}+o(1))n$ contains a 
Hamilton 1-cycle. (The proof in \cite{KKMO} is based on a `hypergraph blow-up 
lemma' due to Keevash~\cite{K}.) This was also proved independently by H\`an and
Schacht~\cite{HS} using a different method. In fact, they showed that if $1 \leq \ell < k/2$,
then any $k$-graph $\Hy$ on $n$ vertices with $(k-\ell)|n$ and 
$\delta(\Hy)\geq(\frac{1}{2(k-\ell)}+o(1))n$ contains a Hamilton $\ell$-cycle. 
They raised the question of determining the correct minimum degree for those 
values of $k$ and $\ell$ not covered by their result or by Theorem~\ref{rrs}. 
Our main result confirms their conjecture and generalises their result. 

\begin{theo}\label{main}
For all $k\geq 3$, $1 \leq \ell \leq k-1$ such that $(k-\ell)\nmid k$ and any 
$\eta > 0$ there exists $n_0$ so that if $n > n_0$ and $(k-\ell)|n$ then any 
$k$-graph $\Hy$ on $n$ vertices with $\delta(\Hy) \geq \left(\frac{1}{\lceil 
\frac{k}{k-\ell} \rceil(k-\ell)}+\eta\right) n$ contains a Hamilton 
$\ell$-cycle.
\end{theo}

This result is best possible up to the error term~$\eta n$, as shown by 
Proposition~\ref{bestpos}. Thus Theorem~\ref{rrs} and Theorem~\ref{main} 
together give asymptotically, for any~$k$ and~$\ell$, the minimum degree 
required to guarantee that a $k$-graph on~$n$ vertices contains a Hamilton 
$\ell$-cycle. {The difference in the minimum degree threshold between the cases $(k-\ell) \mid k$ and $(k-\ell) \nmid k$ is perhaps surprising. For example, if $k=9$ then the minimum degree threshold for an 8-cycle or a 6-cycle is asymptotically $n/2$, whereas for a 7-cycle it is instead $n/10$. This difference is essentially a consequence of the fact that in the $(k-\ell) \mid k$ case every Hamilton $\ell$-cycle contains a perfect matching. The minimum degree threshold for the latter is known to be close to $n/2$ (see Proposition~\ref{rrsbestpos}).} 

Also, less restrictive notions of hypergraph cycles have been considered, e.g. 
in~\cite{B}.

\section{Extremal examples and outline of the proof of Theorem~\ref{main}}

{The next two propositions show that Theorem~\ref{rrs} and Theorem~\ref{main} are each best possible, up to the 
error term~$\eta n$. These constructions are well known, but we include them here for completeness.} By a \emph{perfect matching} in 
a $k$-graph $\Hy$, we mean a set of disjoint edges of $\Hy$ whose union contains 
every vertex of $\Hy$.

\begin{prop} \label{rrsbestpos}
For all~$k\geq 3$, $1 \leq \ell \leq k-1$ and every $n \geq 3k$ such that 
$(k-\ell)|k$ and $k|n$ there exists a $k$-graph $\Hy$ on $n$ vertices with 
$\delta(\Hy) \geq 
\frac{n}{2}-k$ which does not contain a Hamilton $\ell$-cycle.
\end{prop}

\proof Choose $\frac{n}{2} - 1 \leq a \leq \frac{n}{2}+1$ so that $a$ is odd. 
Let $V_1$ and $V_2$ be disjoint sets of size $a$ 
and $n-a$ respectively, and let $\Hy$ be the $k$-graph on vertex set $V = V_1 
\cup V_2$ and with all those $k$-element subsets $S$ of $V$ such that $|S \cap 
V_1|$ 
is even as edges. Then $\delta(\Hy) \geq \min(a, n-a) - k+1 \geq 
\frac{n}{2}-k$. Now, any Hamilton $\ell$-cycle $C$ in $\Hy$ would contain a 
perfect matching, consisting of every $\frac{k}{k-\ell}$th edge of $C$. Every 
edge in this matching would contain an even number of vertices from $V_1$, and 
so $|V_1|$ would be even. Since $|V_1| = a$ is odd, $\Hy$ cannot contain a 
Hamilton $\ell$-cycle.
\endproof

\begin{prop} \label{bestpos}
For all~$k\geq 3$, $1 \leq \ell \leq k-1$ and every $n$ with $(k-\ell)|n$ there 
exists a $k$-graph $\Hy$ on $n$ vertices with $\delta(\Hy) \geq 
\frac{n}{\lceil\frac{k}{k-\ell}\rceil(k-\ell)}-1$ which does not contain a 
Hamilton $\ell$-cycle.
\end{prop}

\proof Let $a:=\lceil\frac{n}{\lceil\frac{k}{k-\ell}\rceil(k-\ell)}\rceil-1$ and 
let $V_1$ and $V_2$ be disjoint sets of size $a$ and $n - a$ respectively. Let 
$\Hy$ be the $k$-graph on vertex set $V = V_1 \cup V_2$ whose edges are all 
those $k$-sets of vertices which contain at least one vertex from~$V_1$. Then 
$\delta(\Hy) = a$. However, an $\ell$-cycle on~$n$ vertices has $n/(k-\ell)$ 
edges and every vertex on such a cycle lies in at most 
$\lceil\frac{k}{k-\ell}\rceil$ edges. Since $\lceil \frac{k}{k-\ell} 
\rceil|V_1|<n/(k-\ell)$, $\Hy$ cannot contain a Hamilton $\ell$-cycle.
\endproof

{A recent construction of Markstr\"om and Ruci\'nski (\cite{MR}) shows that Proposition~\ref{rrsbestpos} still holds if we drop the requirement that $k \mid n$.}

In our proof of Theorem~\ref{main} we construct {a} Hamilton $\ell$-cycle by 
finding several $\ell$-paths and joining them into a spanning $\ell$-cycle. Here 
a $k$-graph $P$ is an $\ell$-\emph{path} if its vertices can be given a linear 
ordering such that every edge of $P$ consists of $k$ consecutive vertices, and 
so that every pair of {consecutive} edges of $P$ (in the natural ordering induced on
the edges) intersect in precisely $\ell$ 
vertices. We say that an enumeration $v_1, v_2, \dots, v_r$ of the vertices of 
$P$ is a \emph{vertex sequence} of $P$ if the edges of $P$ are 
$\{v_{s(k-\ell)+1}, \dots, v_{s(k-\ell)+k}\}$ for each $0 \leq s \leq 
(r-k)/(k-\ell)$. We say that ordered sets $A$ and $B$ are \emph{ordered ends} of 
$P$ if $|A| = |B| = \ell$ and $A$ and $B$ are initial and final segments of a 
vertex sequence of $P$. This allows us to join up $\ell$-paths in the following 
manner. 
Let $P$ and $Q$ be $\ell$-paths, and let $P^{beg}$ and $P^{end}$ be ordered ends 
of $P$, and $Q^{beg}$ and $Q^{end}$ be ordered ends of $Q$. Suppose that 
$P^{end} = Q^{beg}$, and that $V(P) \cap V(Q) = P^{end}$. Then the $k$-graph 
with vertex set $V(P) \cup V(Q)$ and with all the edges of $P$ and of $Q$ is an 
$\ell$-path with ordered ends $P^{beg}$ and $Q^{end}$.

Our proof of Theorem~\ref{main} uses ideas of~\cite{HS}, which in turn were 
based on the `absorbing path' method of~\cite{RRS} and~\cite{RRS2}. Our proof 
contains further developments of the method, which may be of independent 
interest. Roughly 
speaking, the absorbing technique proceeds as follows. We shall prove an 
`absorbing path lemma', which states that in any sufficiently large $k$-graph of 
large minimum degree there exists an $\ell$-path $P$ which can `absorb' any 
small set~$X$ of vertices outside $P$. By this we mean that for any such small 
set $X$ there is another $\ell$-path~$Q$ with the same ordered ends as $P$ and 
with $V(Q) = V(P) \cup X$. Then we can think of replacing~$P$ with $Q$ as 
`absorbing' the vertices of~$X$ into $P$. We shall also prove a `path cover 
lemma', which states that any sufficiently large $k$-graph satisfying the 
minimum degree condition of Theorem~\ref{main} can be almost covered by a 
bounded number of disjoint $\ell$-paths. We can then prove Theorem~\ref{main} 
by combining these lemmas as follows. Firstly, we find in $\Hy$ an absorbing 
$\ell$-path, and then we almost cover the induced $k$-graph on the remaining 
vertices by disjoint $\ell$-paths. We connect up all of these $\ell$-paths to 
form an $\ell$-cycle $C$ which thus contains almost every vertex of $\Hy$. 
Finally, we absorb all vertices of $\Hy$ not contained in $C$ into our absorbing 
path, thereby forming an $\ell$-cycle containing every vertex of $\Hy$.

Beyond these similarities, we have had to make substantial changes to the method 
of H\`an and Schacht. For example, it is simple to `connect up' $\ell$-paths $P$ 
and $Q$ in a $k$-graph $\Hy$ of large minimum degree when $1 \leq \ell < k/2$. 
Indeed, we may add any $k-1-2\ell$ vertices from outside $P$ and $Q$ to the 
ordered ends 
of $P$ and $Q$ to obtain a set $S$ of size $k-1$. Then we can apply the minimum 
degree condition of $\Hy$ to find a vertex $x \in V(\Hy) \sm (V(P) \cup V(Q))$ 
such that $S \cup \{x\}$ is an edge of $\Hy$. Then $P, S \cup \{x\}$ and $Q$ together 
form a single $\ell$-path in~$\Hy$. However, if $\ell \geq k/2$ then things are 
more difficult. So to allow us to connect $\ell$-paths, in 
Section~\ref{sec:diameter} 
we shall use strong hypergraph regularity to prove a `diameter lemma', which 
states that if $1 \leq \ell \leq k-1$ is such that $(k-\ell)\nmid k$, and $A$ 
and $B$ are ordered sets of $\ell$ vertices of a $k$-graph $\Hy$ which has large 
minimum degree, then $\Hy$ contains an $\ell$-path with ordered ends $A$ and $B$ 
with a bounded number of vertices (i.e. the number of vertices depends only on 
$k$). 

In Section~\ref{sec:absorb} we prove our absorbing path lemma. Actually, we will 
not be able to absorb arbitrary sets of vertices, but only `good' $\ell$-sets of 
vertices. We will use strong hypergraph regularity to show that most $\ell$-sets 
of vertices are good, which will be sufficient for our purposes. This weaker 
notion of absorption may be useful for other problems. In Section 
\ref{sec:pathcover} we shall prove the path 
cover lemma. A similar result was already proved in \cite{RRS2}. 
The main difference is that they used weak regularity, whereas we have used 
strong regularity, but this is simply to avoid having to introduce two different 
notions of regularity \textemdash~weak regularity would have sufficed for this 
part of our proof. Finally, in Section \ref{sec:proof} we complete the proof as 
outlined earlier.

\section{Definitions and a preliminary result}

We begin with some notation. By $[r]$ we denote the set of integers from 1 to 
$r$. For a set~$A$, we use $\binom{A}{k}$ to denote the collection of subsets of 
$A$ of size $k$. We write $x = y \pm z$ to denote 
that $y-z \leq x \leq y+z$. By $0 < \alpha \ll 
\beta$ we mean that there exists an increasing function $f: \mathbb{R} 
\rightarrow \mathbb{R}$ such that the following argument is valid for any $0 < 
\alpha \leq f(\beta)$. We write $o(1)$ to denote a function which tends to zero 
as $n$ tends to infinity, holding all other variables involved constant. We 
shall omit floors and ceilings throughout this paper whenever they do not affect 
the argument.

Let $\Hy$ be a $k$-graph on vertex set $V$, with edge set $E$. Then the 
\emph{order} of $\Hy$, denoted~$|\Hy|$, is the number of vertices of $\Hy$ (so 
$|\Hy| = |V|$). For $A 
\subseteq V$, the \emph{neighbourhood} of $A$ is $N_\Hy(A) := \{B \subseteq V: A 
\cup B \in E, A \cap B = \emptyset\}$. The \emph{degree} of $A$, denoted 
$d_\Hy(A)$, is the number of edges of $\Hy$ which contain $A$ as a subset, so 
$d_\Hy(A) = |N_\Hy(A)|$. This is consistent with our previous definition of 
degree for sets of $k-1$ vertices. For any $V' \subseteq V$, the 
\emph{restriction of $\Hy$ to $V'$}, denoted $\Hy[V']$, is the $k$-graph with 
vertex set $V'$ and edges all those edges of $\Hy$ which are subsets of $V'$.

Given two ordered $\ell$-sets of vertices of~$\Hy$, say~$S$ and~$T$, an 
\emph{$\ell$-path from~$S$ to~$T$} in~$\Hy$ is an $\ell$-path in~$\Hy$ which has 
a vertex sequence beginning with the ordered $\ell$-set~$S$ and ending with the 
ordered $\ell$-set~$T$ (i.e. an $\ell$-path with ordered ends $S$ and $T$). 
We say that a $k$-graph $\Hy$ is 
\emph{$s$-partite} if its vertex set $V$ can be partitioned into $s$ vertex 
classes $V_1, \dots, V_s$ such that no edge of $\Hy$ contains more than one 
vertex from any vertex class $V_i$. We denote by $\K[V_1, \dots, V_s]$ the 
\emph{complete $s$-partite $k$-graph} with vertex classes $V_1, \dots, V_s$, 
that is, the $k$-graph with vertex set $V = V_1 \cup \dots \cup V_s$ and edges 
all sets $S \in \binom{V}{k}$ with $|S \cap V_i| \leq 1$ for all $i$.

The following proposition regarding the existence of $\ell$-paths in complete 
$k$-partite 
$k$-graphs will be required in the proof of both the diameter lemma and the 
absorbing path lemma.

\begin{prop} \label{pathsincompletegraph}
Suppose that $k \geq 3$, and that $1 \leq \ell \leq k-1$ is such that $(k-\ell) 
\nmid k$. Let $V$ be a set 
of vertices partitioned into $k$ vertex classes $V_1, ..., V_k$, with $|V_i| = 
k\ell (k-\ell)+1$ for each~$i$, and let $P^{beg}$ and $P^{end}$ be disjoint 
ordered 
sets of $\ell$~vertices from $V$ such that $|P^{beg} \cap V_i| \leq 1$ and 
$|P^{end} \cap 
V_i| \leq 1$ for each $1 \leq i \leq k$. Then $\K[V_1, ..., V_k]$ contains an 
$\ell$-path $P$ from $P^{beg}$ to $P^{end}$ containing every vertex of $V$ (so 
$|V(P)| = k^2\ell(k-\ell)+k$).
\end{prop}

\proof To prove this result, we consider strings (finite sequences of 
characters) on character set $[k]$. We denote the $i$th character of a string 
$S$ by $S_i$. By an ordering of $[k]$ we mean a string of length $k$ which 
contains each character precisely once. Let $A$ and $B$ be orderings of~$[k]$. 
We say that $A$ and $B$ are adjacent if we can obtain $B$ from $A$ by swapping a 
single pair of adjacent characters in $A$. So for example, $12345$ is adjacent 
to $12435$.

Suppose that $A$ and $B$ are adjacent orderings of $[k]$, and let $i$ and $i+1$ 
be the positions in $A$ of the characters swapped to obtain $B$ from $A$ (so $1 
\leq i \leq k-1$). Since $(k-\ell) \nmid k$ we may choose $p \in \{1,2\}$ such 
that $(k-\ell) \nmid ((p-1)k+i)$. Then define the string $S(A,B)$ to consist of 
$p$ consecutive copies of $A$ followed by $(k-\ell+1)-p$ copies
of $B$. Then 
$S(A,B)$ has length $(k-\ell+1)k$ and the property that $S(A,B)$ starts with $A$ 
and ends with $B$. Note that the only consecutive subsequence of $S$ of length 
$k$ which contains some character more than once is $S' = 
S(A,B)_{(p-1)k+i+1}\dots S(A,B)_{pk+i}$. In other words, $S'$ contains the final 
$k-i$ characters of $A$ and the first $i$ characters of $B$, and the first and final 
character of $S'$ is $A_{i+1}$. Therefore, as $(k-\ell) \nmid ((p-1)k+i)$, we 
know that no character appears twice in $S(A,B)_{r(k-\ell)+1}, \dots, 
S(A,B)_{r(k-\ell)+k}$ for any $0 \leq r \leq k$. Furthermore, $S(A,B)$ contains 
the same number of copies of each character.

Now, choose a string $C$ to be any ordering of $[k]$ such that for $1 \leq i 
\leq \ell$, the $i$th vertex of the ordered set $P^{beg}$ lies in vertex class 
$V_{C_i}$. Define a string $D$ to be an ordering of $[k]$ such that for 
$1 \leq i \leq \ell$, the $i$th vertex of the ordered set $P^{end}$ lies in 
vertex class $V_{D_{i+k-\ell}}$, and the characters $D_i$ for $1 \leq i \leq 
k-\ell$ appear in the same order as they do in $C$. Then we may transform $C$ into $D$ 
through at most $k\ell$ swaps of pairs of consecutive vertices. So we may choose 
$A^0, \dots, A^{k\ell}$ to be orderings of $[k]$ such that $A^0 = C$, $A^{k\ell} 
= D$, and for any $0 \leq i \leq k\ell-1$, $A^i$ and $A^{i+1}$ are either 
adjacent or identical.

Then for each $0 \leq i \leq k\ell-1$ we may choose a string $S^i$ of length 
$k(k-\ell+1)$ such that~$S^i$ starts with $A^i$ and ends 
with $A^{i+1}$, each character appears an equal number of times in~$S^i$ and for 
each $0 \leq r \leq k$ no character appears more than once in 
$S^i_{r(k-\ell)+1}, 
\dots, S^i_{r(k-\ell)+k}$. Indeed, if $A^i$ and $A^{i+1}$ are adjacent, take 
$S^i$ to be $S(A^i, A^{i+1})$, and if $A^i = A^{i+1}$, take $S^i$ to be the 
string consisting of $k-\ell+1$ 
consecutive copies of $A^i$. For each $0 \leq i \leq k\ell-2$, let $T^i$ be the 
string obtained by deleting the final $k$ characters of $S^i$, and let 
$T^{k\ell-1} = S^{k\ell-1}$. Let $S$ be the string formed by concatenating $T^0, 
\dots, T^{k\ell-1}$. Then $S$ starts with $C$ and ends with~$D$ and has the 
property that no character appears twice in $S_{r(k-\ell)+1}, \dots, 
S_{r(k-\ell)+k}$ for any $0 \leq r \leq k^2\ell$. Also $|S| = 
k^2\ell(k-\ell)+k$, 
and so since $S$ contains each character the same number of times, each 
character 
appears $k\ell(k-\ell)+1$ times in $S$.

We can now construct the vertex sequence of our desired $\ell$-path $P$. To do 
so, let $P$ have vertex sequence beginning with $P^{beg}$ and ending with 
$P^{end}$. In between, let the $i$th vertex of $P$ be chosen from $V_{S_i}$, and 
make these choices without choosing the same vertex twice. Then $P$ contains all 
$k\ell(k-\ell)+1$ vertices from each vertex class and is an $\ell$-path. Indeed, 
the edges of an $\ell$-path $P$ consist of the vertices in positions 
$r(k-\ell)+1, \dots, r(k-\ell)+k$ for $0 \leq r \leq |E(P)| -1$. So by 
construction these vertices are from different vertex classes, and so form an 
edge in $\K[V_1, \dots, V_k]$.\endproof

{Note that Proposition~\ref{pathsincompletegraph} would not hold if instead we had $(k-\ell) \mid k$, as it would not be possible to choose $p$ as in the proof.}

\section{The regularity lemma for $k$-graphs} \label{regularity}

\subsection{Regular complexes}
Before we can state the regularity lemma, we first have to say what we mean by a 
regular or `quasi-random' hypergraph and, more generally, by a regular 
complex. A \emph{hypergraph} $\Hy$ consists of a vertex set $V(\Hy)$ and an 
edge set $E(\Hy)$, where every edge $e \in E(\Hy)$ is a non-empty subset of $V(\Hy)$. So a 
$k$-graph as defined earlier is a hypergraph in which every edge has size $k$. A 
hypergraph~$\Hy$ is a \emph{complex} if whenever~$e\in E(\Hy)$ and $e'$ is a 
non-empty subset of $e$ we have that $e'\in E(\Hy)$. All the complexes 
considered 
in this paper have the property that every vertex forms an edge. A complex~$\Hy$ 
is a \emph{$k$-complex} if every edge of $\Hy$ consists of at most $k$ vertices. 
The edges of size~$i$ are called \emph{$i$-edges} of~$\Hy$. We write 
$|\Hy|:=|V(\Hy)|$ for the \emph{order} of~$\Hy$. Given a $k$-complex $\Hy$, for 
each 
$i=1,\ldots,k$ we denote by~$\Hy_i$ the \emph{underlying $i$-graph of~$\Hy$}. So 
the vertices of~$\Hy_i$ are those of~$\Hy$ and the edges of~$\Hy_i$ are the 
$i$-edges of~$\Hy$.

Note that a $k$-graph $\Hy$ can be turned into a $k$-complex, which we denote by 
$\Hy^\leq$, by making every edge into a complete $i$-graph~$K_k^{(i)}$, for each 
$1 \leq i\leq k$. (In a more general $k$-complex we may have $i$-edges which do 
not lie within an $(i+1)$-edge.) Given $k\le s$, a \emph{$(k,s)$-complex}~$\Hy$ 
is an $s$-partite $k$-complex, by which we mean that the vertex set of~$\Hy$ can 
be partitioned into sets $V_1,\dots,V_s$ such that every edge of~$\Hy$ meets 
each~$V_i$ in at most one vertex.

Given $i \geq 2$, an $i$-partite $i$-graph~$\Hy_i$, and an $i$-partite 
$(i-1)$-graph~$\Hy_{i-1}$ on the same vertex set, we write $\K_i(\Hy_{i-1})$ for 
the set of $i$-sets of vertices which form a copy of the complete 
$(i-1)$-graph~$K_i^{(i-1)}$ on~$i$ vertices in~$\Hy_{i-1}$. We define the 
\emph{density of $\Hy_i$ with respect to $\Hy_{i-1}$} to be
$$ d(\Hy_i|\Hy_{i-1}):= \frac{|\K_i(\Hy_{i-1})\cap E(\Hy_i)|}{|\K_i(\Hy_{i-1})|} 
$$
if $|\K_i(\Hy_{i-1})|>0$, and $d(\Hy_i|\Hy_{i-1}):=0$ otherwise. More generally, 
if 
$\mathbf{Q}:=(Q(1),Q(2),\ldots,Q(r))$ is a collection of~$r$ subhypergraphs 
of~$\Hy_{i-1}$, we define $\K_i(\mathbf{Q}):=\bigcup_{j=1}^r \K_i(Q(j))$ and
$$ d(\Hy_i |\mathbf{Q}):= \frac{|\K_i(\mathbf{Q})\cap 
E(\Hy_i)|}{|\K_i(\mathbf{Q})|} $$ 
if $|\K_i(\mathbf{Q})|>0$, and $d(\Hy_i |\mathbf{Q}):=0$ otherwise. 

We say that $\Hy_i$ is \emph{$(d_i,\delta,r)$-regular with respect to~$\Hy_{i-1}$} 
if every $r$-tuple $\mathbf{Q}$ with 
$|\K_i(\mathbf{Q})|>\delta|\K_i(\Hy_{i-1})|$ 
satisfies $d(\Hy_i|\mathbf{Q}) = d_i \pm \delta. $
Instead of $(d_i,\delta,1)$-regularity we sometimes refer to 
\emph{$(d_i,\delta)$-regularity}.

Given $3\le k\le s$ and a $(k,s)$-complex~$\Hy$, we say that~$\Hy$ is 
\emph{$(d_k,\ldots,d_2,\delta_k,\delta,r)$-regular} if the following conditions hold:
\begin{itemize} 
\item For every $i=2,\dots,k-1$ and for every $i$-tuple~$K$ of vertex classes 
either $\Hy_i[K]$ is $(d_i,\delta)$-regular with respect to~$\Hy_{i-1}[K]$ or 
$d(\Hy_i[K] | \Hy_{i-1}[K])=0$. 
\item For every $k$-tuple~$K$ of vertex classes either $\Hy_k[K]$ is 
$(d_k,\delta_k,r)$-regular with respect to~$\Hy_{k-1}[K]$ or $d(\Hy_k[K] | 
\Hy_{k-1}[K])=0$.
\end{itemize}
Here we write $\Hy_i[K]$ for the restriction of $\Hy_i$ to the union of all 
vertex 
classes in~$K$. We sometimes denote $(d_k,\ldots,d_2)$ by $\mathbf{d}$ and refer 
to $(\mathbf{d},\delta_k,\delta,r)$-regularity.

We will need the following lemma which states that the restriction of regular 
complexes to a sufficiently large set of vertices is still regular.

\begin{lemma}\label{restriction}
Let $k,s,r,m$ be positive integers and $\alpha,d_2,\dots,d_k,\delta,\delta_k$ be 
positive constants such that
$$1/m\ll 1/r,\delta\le \min\{\delta_k,d_2,\dots,d_{k-1}\}\le \delta_k\ll \alpha\ll 
d_k,1/s.$$
Let $\Hy$ be a $(\mathbf{d}, \delta_k, \delta, r)$-regular $(k,s)$-complex with 
vertex classes $V_1,\dots,V_s$ of size~$m$. 
For each~$i$ let $V'_i\subseteq V_i$ be a set of size at least $\alpha m$. Then 
the restriction $\Hy' = \Hy[V'_1\cup \dots\cup V'_s]$ of~$\Hy$ to $V'_1\cup 
\dots\cup V'_s$ is $(\mathbf{d},\sqrt{\delta_k},\sqrt{\delta},r)$-regular.
\end{lemma}

It is easy to prove Lemma~\ref{restriction} by induction on~$i$ (where $2\le i\le k$
is as in the definition of a regular complex). In the 
induction step, use the dense hypergraph counting lemma (Corollary~6.11 in~\cite{KRS})
to show that $\sqrt {\delta} |\K_i(\Hy'_{i-1})| \geq \delta |\K_i(\Hy_{i-1})|$ and
likewise when $i=k$.

\subsection{Statement of the regularity lemma}
In this section we state the version of the regularity lemma for $k$-graphs due 
to R\"odl and Schacht~\cite{RS}, which we will use several times in our proof. 
To prepare for this we will first need some more notation. Suppose that $V$ is a 
finite set of vertices and $\Part^{(1)}$ is a partition of $V$ into sets 
$V_1,\dots,V_{a_1}$, which will be called \emph{clusters}. Given $k\ge 3$ and
any $j \in [k]$, we denote by $\cross_j=\cross_j(\Part^{(1)})$ the set of all those 
$j$-subsets of~$V$ that meet each~$V_i$ in at most~1 vertex. For every set 
$A\subseteq [a_1]$ with $2\le |A| \leq k-1$ we write $\cross_A$ for all those 
$|A|$-subsets of~$V$ that meet each~$V_i$ with $i\in A$. Let $\Part_A$ be a 
partition of $\cross_A$. We refer to the partition classes of $\Part_A$ as 
\emph{cells}. For each $i=2,\dots,k-1$ let $\Part^{(i)}$ be the union of all the 
$\Part_A$ with $|A|=i$. So $\Part^{(i)}$ is a partition of $\cross_i$. 

$\Part(k-1)=\{\Part^{(1)},\ldots, \Part^{(k-1)} \}$ is a \emph{family of 
partitions on~$V$} if the following condition holds. Recall that $a_1$ denotes the number 
of clusters in $\Part^{(1)}$. Consider any $B\subseteq A\subseteq [a_1]$ such 
that $2\leq |B|<|A| \leq k-1$ and suppose that $S,T\in \cross_A$ lie in the same 
cell of $\Part_A$. Let $S_B:=S\cap\bigcup_{i\in B} V_i$ and define $T_B$ 
similarly. Then $S_B$ and $T_B$ lie in the same cell of~$\Part_B$.

To illustrate this condition, suppose that $k=4$ and $A=[3]$. Then 
$\Part_{\{1,2\}}$, $\Part_{\{2,3\}}$ and $\Part_{\{1,3\}}$ partition the edges 
of the $3$ complete bipartite graphs induced by the pairs $V_1V_2$, $V_2V_3$ and 
$V_1V_3$. These the partitions together naturally induce a partition $\Q$ of the 
set of triples induced by $V_1,V_2$ and $V_3$. The above condition says that 
$\Part_{\{1,2,3\}}$ must be a refinement of~$\Q$.

Given $1 \leq i \le j\le k$ with $i<k$, $J\in \cross_j$ and an $i$-set $Q\subseteq J$, we 
write $C_Q$ for the set of all those $i$-sets in $\cross_{i}$ that lie in the 
same cell of $\Part^{(i)}$ as~$Q$. (In particular, if $i=1$ then~$C_Q$
is the cluster containing the unique element in~$Q$.) The
\emph{polyad $\hat{P}^{(i)}(J)$ of $J$} 
is defined by $\hat{P}^{(i)}(J) := \bigcup_{Q} C_Q$, where the union is over all 
$i$-subsets~$Q$ of~$J$. So we can view $\hat{P}^{(i)}(J)$ as an $j$-partite 
$i$-graph (whose vertex classes are the clusters intersecting~$J$). We let 
$\hat{\Part}^{(j-1)}$ be the set consisting of all the $\hat{P}^{(j-1)}(J)$ for 
all $J\in \cross_j$. So for each $K\in \cross_k$ we can view 
$\bigcup_{i=1}^{k-1}\hat{P}^{(i)}(K)$ as a $(k-1,k)$-complex.

We say that $\Part = \Part(k-1)$ is \emph{$(\eta, \delta,t)$-equitable} if
\begin{itemize}
\item there exists ${\mathbf{d}}=( d_{k-1}, \dots, d_2)$ such that $d_i \geq 
1/t$ and $1/d_i \in \mathbb{N}$ for all $i=2, \dots, k-1$,
\item $\Part^{(1)}$ is a partition of~$V$ into~$a_1$ clusters of equal size, 
where $1/\eta \leq a_1 \leq t$,
\item for all $i = 2, \dots, k-1$, $\Part^{(i)}$ is a partition of $\cross_i$ 
into at most $t$ cells,
\item for every $K \in \cross_k$, the $(k-1,k)$-complex 
$\bigcup_{i=1}^{k-1}\hat{P}^{(i)}(K)$ is $({\bf d},\delta,\delta,1)$-regular.
\end{itemize}
Note that the final condition implies that for all $i=2, \dots, k-1$ the cells 
of $\Part^{(i)}$ have almost equal size.

Let $\delta_k>0$ and~$r\in\mathbb{N}$. Suppose that~$\Hy$ is a $k$-graph on~$V$ 
and $\Part = \Part (k-1)$ is a family of partitions on~$V$. Given a polyad
$\hat{P}^{(k-1)}\in \hat{\Part}^{(k-1)}$, we say that $\Hy$ is 
\emph{$(\delta_k,r)$-regular with respect to $\hat{P}^{(k-1)}$} if~$\Hy$ is 
$(d,\delta_k,r)$-regular with respect to $\hat{P}^{(k-1)}$ for some~$d$. We say 
that~$\Hy$ is \emph{$(\delta_k,r)$-regular with respect to~$\Part$} if
$$\left| \bigcup \{ \K_k(\hat{P}^{(k-1)}): \; \mbox{ $\Hy$ is not 
$(\delta_k,r)$-regular with respect to $\hat{P}^{(k-1)} \in 
\hat{\Part}^{(k-1)}$}  \} \right| \leq \delta_k |V|^k.$$
This means that not much more than a $\delta_k$-fraction of the $k$-subsets 
of~$V$ form a~$K_k^{(k-1)}$ that lies within a polyad with respect to 
which~$\Hy$ 
is not regular.

Now we are ready to state the regularity lemma.

\begin{theo}[R\"odl and Schacht~\cite{RS}, Theorem 17]\label{reglemma}
Let $k\geq 3$ be a fixed integer. For all positive constants~$\eta$ 
and~$\delta_k$ and all functions $r: {\mathbb N}\rightarrow \mathbb{N}$ 
and $\delta: \mathbb{N} \rightarrow (0,1]$, there are integers~$t$ 
and~$n_0$ such that the following holds for all $n \ge n_0$ which are divisible 
by~$t!$. Suppose that~$\Hy$ is a $k$-graph of order~$n$. Then there exists a 
family of partitions $\Part = \Part(k-1)$ of the vertex set~$V$ of~$\Hy$ such 
that
\begin{enumerate}
\item $\Part$ is $(\eta, \delta(t),t)$-equitable and
\item $\Hy$ is $(\delta_k,r(t))$-regular with respect to~$\Part$.
\end{enumerate}
\end{theo}

Similar results were proved earlier by R\"odl and Skokan~\cite{RSk} and 
Gowers~\cite{G}.
Note that the constants in Theorem~\ref{reglemma} can be chosen so that they 
satisfy the following hierarchy:
\begin{equation*} \label{reghierarchy}
\frac{1}{n_0} \ll \frac{1}{r}=\frac{1}{r(t)},\delta=\delta(t)
\ll \min\{\delta_k,1/t\} \ll \eta.
\end{equation*}

\subsection{The reduced $k$-graph}\label{sec:rhgraph}

To prove the absorbing lemma and the path cover lemma, we will use the so-called 
reduced $k$-graph. Suppose that we have constants
\begin{equation*}
\frac{1}{n_0} \ll \frac{1}{r},\delta \ll \min\{\delta_k,1/t\} \leq \delta_k, 
\eta \ll d \ll \theta \ll \mu, 1/k.
\end{equation*}
and a $k$-graph $\Hy$ on $V$ of order $n \geq n_0$ with $\delta(\Hy) \geq (\mu + 
\theta) n$. We may apply the regularity lemma to $\Hy$ to obtain a family of 
partitions $\Part = \{ \Part^{(1)},\ldots, \Part^{(k-1)} \}$ of~$V$. Then the 
\emph{reduced $k$-graph} $\R=\R(\Hy,\Part)$ is the $k$-graph whose vertices are 
the clusters of $\Hy$, i.e.~the parts of~$\Part^{(1)}$. A $k$-tuple of clusters 
forms an edge of $\R$ if there is some polyad~$\hat{P}^{(k-1)}$ induced on these 
$k$ clusters such that $\Hy$ is $(d', \delta_k, r)$-regular with respect to 
$\hat{P}^{(k-1)}$ for some $d' \geq d$. To make use of the reduced $k$-graph, we 
shall need to show that it almost inherits the minimum degree condition from 
$\Hy$.

\begin{lemma}\label{minreduced}
All but at most $\theta |\R|^{k-1}$ sets $S \in \binom{V(\R)}{k-1}$ satisfy 
$d_\R(S) \geq \mu |\R|$.
\end{lemma}

Similar results have been proved in previous papers on hypergraph Hamilton 
cycles, but we include the short proof for completeness, in which we will need 
the following lemma. We say 
that an edge~$e$ of $\Hy$ is \emph{useful} if it lies in $\K_k(\hat{P}^{(k-1)})$ for 
some $\hat{P}^{(k-1)} \in \hat{\Pa}^{(k-1)}$ such 
that $\Hy$ is $(d', \delta_k, r)$-regular with respect to $\hat{P}^{(k-1)}$ for 
some $d' \geq d$. Note that if $e$ lies in $\K_k(\hat{P}^{(k-1)})$ then
$\hat{P}^{(k-1)}=\hat{P}^{(k-1)}(e)$ is the polyad of~$e$. Moreover, if $e$
is a useful edge of $\Hy$, and $V_{i_1}, \dots, V_{i_k}$ are the clusters
containing the vertices of $e$, then these $k$ clusters will form an edge of $\R$. 

\begin{lemma} \label{mostedgesgood}
At most $2d n^k$ edges of $\Hy$ are not useful.
\end{lemma}

\proof
There are three reasons why an edge of $\Hy$ may not be useful. Firstly, it may lie 
in $\binom{V}{k} \setminus \cross_k$. Since $\Part^{(1)}$ partitions $V$ into $a_1$ clusters 
of equal size, there are at most $ \frac{n}{a_1} n^{k-1} \leq \eta n^k$ edges of 
this type. Secondly, the edge may lie in a polyad $\hat{P}^{(k-1)} \in 
\hat{\Part}^{(k-1)}$ such that $|E(\Hy) \cap \K_k(\hat{P}^{(k-1)})| \leq 
d|\K_k(\hat{P}^{(k-1)})|$. There are at most $dn^k$ edges of this type. Finally, 
the edge may lie in a polyad $\hat{P}^{(k-1)} \in \hat{\Part}^{(k-1)}$ such that 
$\Hy$ is not $(\delta_k,r)$-regular with respect to $\hat{P}^{(k-1)}$.
Since $\Hy$ is $(\delta_k, r)$-regular with respect to 
$\Part$, there are at most $\delta_k n^k$ edges of this type. So altogether, at 
most $(\delta_k+d+\eta)n^k \leq 2dn^k$ edges of $\Hy$ are not useful.
\endproof

\removelastskip\penalty55\medskip\noindent{\bf Proof of Lemma~\ref{minreduced}.}
Let $m = |V_1| = \dots = |V_{a_1}|$ be the size of the clusters. We say that a 
$(k-1)$-tuple of clusters of $\Hy$ is \emph{poor} if there are at least $\theta m^{k-1} 
n$ edges of $\Hy$ which intersect each of the $k-1$ clusters in precisely one 
vertex and which are 
not useful. Then it follows from Lemma~\ref{mostedgesgood} that at most $\theta 
|\R|^{k-1}$ such $(k-1)$-tuples are poor. So it remains to show that any
$(k-1)$-tuple which is not poor has many neighbours in $\R$. But if $V_{i_1}, \dots, V_{i_{k-1}}$ 
is a $(k-1)$-tuple which is not poor, then there are at least
$m^{k-1}\delta(\Hy)-\theta m^{k-1}n\ge \mu m^{k-1}n$ useful edges of 
$\Hy$ which intersect each of $V_{i_1}, \dots, V_{i_{k-1}}$ in precisely one 
vertex. For any other cluster~$V_j$ at most $m^k$ edges of $\Hy$ intersect each 
of $V_{i_1}, \dots, V_{i_{k-1}}, V_j$ in precisely one vertex, and so there are 
at least $\mu 
n/m = \mu |\R|$ choices of $V_j$ such that there is at least one such 
useful edge. 
This useful edge indicates the existence of a polyad satisfying the conditions 
of an edge in the reduced $k$-graph~$\R$.
\endproof

\subsection{The embedding and the extension lemmas.}

In our proof we will also use an embedding lemma, which guarantees the 
existence of a copy of a complex~$\G$ of 
bounded maximum degree inside a 
suitable regular complex~$\Hy$, where the order of~$\G$ is allowed to be linear 
in the order of~$\Hy$. In order to state this lemma, we need some more 
definitions.

The degree of a vertex~$x$ in a complex~$\G$ is the number of edges 
containing~$x$. The \emph{maximum vertex degree} of~$\G$ is the largest degree 
of a vertex of~$\G$.

Suppose that~$\Hy$ is a $(k,s)$-complex with vertex classes $V_1,\dots,V_s$, 
which all have size~$m$. Suppose also that~$\G$ is a $(k,s)$-complex with 
vertex classes $X_1,\dots,X_s$ of size at most~$m$. We say that~$\Hy$ 
\emph{respects the partition of~$\G$} if whenever~$\G$ contains an $i$-edge 
with vertices in $X_{j_1},\ldots,X_{j_i}$, then there is an $i$-edge of~$\Hy$ 
with vertices in $V_{j_1},\ldots,V_{j_i}$. On the other hand, we say that a 
labelled copy of~$\G$ in~$\Hy$ is \emph{partition-respecting} if for each 
$i=1,\dots,s$ the vertices corresponding to those in~$X_i$ lie within~$V_i$. 

\begin{lemma}[Embedding lemma, \cite{CFKO}, Theorem 3]\label{emblemma}
Let $\Delta,k,s,r,m_0$ be positive integers and let 
$c,d_2,\ldots,d_k,\delta,\delta_k$ be positive constants such that $1/d_i 
\in\mathbb{N}$ for all $i < k$, 
$$1/m_0\ll 1/r, \delta \ll\min\{\delta_k,d_2,\ldots,d_{k-1}\}\le\delta_k\ll 
d_k,1/\Delta,1/s$$
and
$$c\ll d_2,\ldots,d_k.$$
Then the following holds for all integers $m\ge m_0$. Suppose that $\G$ is a 
$(k,s)$-complex of maximum vertex degree at most~$\Delta$ with vertex classes 
$X_1,\dots,X_s$ such that $|X_i|\le cm$ for all $i=1,\dots,s$. Suppose also that 
$\Hy$ is a $(\mathbf{d},\delta_k,\delta,r)$-regular $(k,s)$-complex with vertex 
classes $V_1,\dots,V_s$, all of size $m$, which respects the partition of~$\G$.
Then $\Hy$ contains a labelled partition-respecting copy of $\G$.
\end{lemma}

We will also use the following weak version of a lemma 
from~\cite{CFKO}. Roughly speaking, it states that if $\G$ is an induced 
subcomplex of $\G'$, and $\Hy$ is suitably regular, then almost all copies
of $\G$ in $\Hy$ can be extended to a large number of copies of $\G'$ in $\Hy$. We 
write $|\G|_\Hy$ for the number of labelled partition-respecting copies 
of~$\G$ in~$\Hy$.

\begin{lemma}[Extension lemma, \cite{CFKO},  Lemma 5]\label{extensions_count}
Let $k,s,r,b',b'',m_0$ be positive integers, where $b'<b''$, and let $c, 
\beta,d_2,\ldots,d_k,\delta,\delta_k$ be positive constants such that $1/d_i 
\in\mathbb{N}$ for all $i <k$ and 
$$1/m_0\ll 1/r, \delta \ll c \ll \min\{\delta_k,d_2,\ldots,d_{k-1}\} \le 
\delta_k \ll \beta,d_k,1/s,1/b''.$$
Then the following holds for all integers $m\ge m_0$. Suppose that $\G'$ is a 
$(k,s)$-complex on $b''$ vertices with vertex classes $X_1,\dots,X_s$ and 
let~$\G$ be an induced subcomplex of~$\G'$ on~$b'$ vertices. Suppose also that 
$\Hy$ is a $(\mathbf{d},\delta_k,\delta,r)$-regular $(k,s)$-complex with vertex 
classes $V_1,\dots,V_s$, all of size $m$, which respects the partition of 
$\G'$. Then all but at most $\beta |\G|_\Hy$ labelled partition-respecting 
copies of~$\G$ in~$\Hy$ are extendible to at least $cm^{b''-b'}$ labelled 
partition-respecting copies of~$\G'$ in~$\Hy$.
\end{lemma}

The proofs of Lemmas~\ref{emblemma} and~\ref{extensions_count} rely on the 
hypergraph counting lemma (Theorem~9 in~\cite{RS2}). In particular, the extension lemma is a 
straightforward consequence of the counting lemma.
Actually both the embedding lemma and the extension lemma involved the 
additional condition that $1/{d_k} \in \mathbb{N}$. However, this can easily be 
achieved by working with a subcomplex~$\Hy'$ of $\Hy$ which is 
$(d'',d_{k-1},\dots,d_2,\delta_k,\delta,r)$-regular with respect to
$\hat{P}^{(k-1)}$ for some $d''\gg \delta_k$ with $1/d''\in \mathbb{N}$.
The existence of such a $\Hy'$ follows 
immediately from the slicing lemma (\cite{RS}, Proposition 22), which is proved 
using a simple application of a Chernoff bound.

Now suppose that we have applied the regularity lemma (Theorem~\ref{reglemma}) 
to a $k$-graph~$\Hy$ to obtain a reduced $k$-graph~$\R$. An edge~$e$ of~$\R$ 
indicates that we can apply the embedding lemma or the extension lemma to the 
subcomplex of~$\Hy$ whose vertex classes are the 
clusters $V_1,\dots,V_k$ corresponding to the vertices of~$e$. More precisely, 
since $e$ 
is an edge of~$\R$, there is some polyad $\hat{P}^{(k-1)}=\hat{P}^{(k-1)}(K)$ 
(where $K\in {\rm Cross}_k$) induced by $V_1,\dots,V_k$ such that $\Hy$ is 
$(d',\delta_k,r)$-regular with respect to $\hat{P}^{(k-1)}$ for some $d'\ge d$. 
Let $\Hy^*$ be the $(k,k)$-complex obtained from the $(k-1,k)$ complex 
$\bigcup_{i=1}^{k-1}\hat{P}^{(i)}(K)$ by adding~$E(\Hy)\cap \K(\hat{P}^{(k-1)})$ 
as the `$k$th level'. Then $\Hy^*$ is a $(\mathbf{d},\delta_k,\delta,r)$-regular 
subcomplex of $\Hy$, where ${\bf d} =(d',d_{k-1}, \dots, d_2)$, and $(d_{k-1}, 
\dots, d_2)$ 
is as in the definition of a $(\eta, \delta, t)$-equitable family of partitions. 
Also~$\Hy^*$ satisfies the 
conditions of the embedding (or extension) lemma. So in particular, the 
embedding lemma implies that if $m:=|V_1|$ and $\G$ is a $k$-partite $k$-graph 
of bounded maximum vertex degree whose vertex classes have size at most $cm$, 
then $\Hy$ contains a copy of~$\G$.

\section{Diameter Lemma} \label{sec:diameter}

In this section, we shall prove a diameter lemma, which will state that any 
sufficiently large $k$-graph of large minimum degree {has small diameter in the sense that we can find an $\ell$-path 
from any ordered $\ell$-set of vertices to any other ordered $\ell$-set of 
vertices. A similar assertion for the case $\ell = k-1$, (called the `Connecting Lemma'), was proved in \cite{RRS2}. The proof is quite different from ours.} To prove the diameter lemma, we shall first consider a $k$-graph $\W(k, \ell)$, for 
which a similar statement is easier to prove 
(Proposition~\ref{connectedclusters}). For $k/2 \leq \ell \leq k-2$, the 
$k$-graph $\W(k, \ell)$ has $4\ell - k +2$ vertices in three disjoint sets $X, 
Y$~and~$Z$, where $X = \{x_1, \dots, x_\ell\}, Y = \{y_1, \dots, y_\ell\}$ and 
$Z = \{z_1, \dots, z_{2\ell - k +2}\}$. $\W(k, \ell)$ has $2\ell -k +2$ edges, 
where for $1 \leq i \leq 2\ell-k+2$ the $i$th edge of $\W(k, \ell)$ is $\{x_1, 
\dots, x_{\ell+1-i}\} \cup \{y_1, \dots, y_{k-2-\ell+i}\} \cup \{z_i\}$. So each 
edge of $\W(k,\ell)$ intersects the following edge in precisely $k-2$ vertices. 
We shall sometimes view $\W(k, \ell)$ as a $(4\ell - k +2)$-partite $k$-graph 
with a single vertex in each vertex class, and consider the $(k, 4\ell - k 
+2)$-complex $\W(k, \ell)^\leq$. We refer to the ordered sets $X$ and~$Y$ as 
the ordered ends of $\W(k, \ell)$.

The next proposition states that for most pairs of sets $S$ and $T$ of $\ell$ vertices 
in a $k$-graph $\Hy$ of large minimum degree, $\Hy$ contains many copies of 
$\W(k, \ell)$ with $S$ and $T$ as ordered ends.

\begin{prop} \label{connectedclusters}
Suppose that $k \geq 3$, that $k/2 \leq \ell \leq k-2$ and that $1/n \ll \gamma 
\ll \beta \ll \mu, 1/k$. Let~$\Hy$ be a $k$-graph on~$n$ vertices such that 
$d(S) \geq \mu n$ for all but at most~$\gamma n^{k-1}$ sets~$S \in 
\binom{V(\Hy)}{k-1}$. Then for all but at most~$\beta n^{2\ell}$ pairs~$S,T$ of 
ordered $\ell$-sets of vertices of $\Hy$ there are at least $\beta n^{2\ell-k+2}$ copies 
of~$\W(k, \ell)$ in~$\Hy$ with ordered ends~$S$ and~$T$.
\end{prop}

\proof We refer to the at most $\gamma n^{k-1}$ sets $S$ of $k-1$ vertices in 
$\Hy$ which do not satisfy $d(S) \geq \mu n$ as unfriendly $(k-1)$-sets. We say 
that a 
pair of $\ell$-sets $S$ and $T$ is unfriendly if there exist $S' \subseteq S$, 
$T' \subseteq T$ such that $S' \cup T'$ is a unfriendly $(k-1)$-set. Then for any 
unfriendly $(k-1)$-set $B$, there are at most $2^{k-1}n^{2\ell-k+1}$ pairs of $\ell$-sets 
$S$ and $T$ with $S' \cup T' = B$ for some $S' \subseteq S$ and $T' \subseteq 
T$, and so since there are at most $\gamma n^{k-1}$ unfriendly $(k-1)$-sets, and 
$\gamma \ll \beta\ll 1/k$, we know that there are at most $\beta n^{2\ell}$ unfriendly 
pairs of $\ell$-sets. 

To complete the proof, it is sufficient to show that if the pair $S, T$ of 
ordered $\ell$-sets is not unfriendly, then $\Hy$ contains at least $\beta 
n^{2\ell-k+2}$ copies of~$\W(k,\ell)$ 
with ordered ends~$S$ and~$T$. Let $S = \{x_1, \dots, x_\ell\}$, and let $T = \{y_1, 
\dots, y_\ell\}$. For each $1 \leq i \leq 2k-\ell+2$ we choose a vertex $z_i$ 
such that $z_i \notin S \cup T$, $z_i \neq z_j$ for any $j < i$, and such that 
$\{x_1, \dots, x_{\ell+1-i}, y_1, \dots, y_{k-2-\ell+i}, z_i\}$ is an edge of 
$\Hy$. This is possible for each $i$ as we know that $S,T$ is not a unfriendly 
pair, 
and so $d(\{x_1, \dots, x_{\ell+1-i}, y_1, \dots, y_{k-2-\ell+i}\}) \geq \mu n$, 
and hence there are at least $\mu n - (4\ell-k+2)$ vertices to choose from. Then 
$S, T$ and the chosen vertices $z_i$ together form a copy of $\W(k, \ell)$ in 
$\Hy$ with ordered ends $S$ and~$T$. Since $\beta \ll \mu$, by counting the choices we 
could have made for the $z_i$ we find that $\Hy$ contains at least $\beta n^{2\ell-k+2}$ 
copies of $\W(k, \ell)$ with ordered ends $S$ and $T$. \endproof

The following proposition relates the $k$-graph $\W(k, \ell)$ to a $k$-graph 
$\Pa(k, \ell)$ which consists of several $\ell$-paths from one ordered 
$\ell$-set to another. We say that $\ell$-paths $P$ and $Q$ with ordered ends 
$P^{beg}$, $P^{end}$, $Q^{beg}$ and $Q^{end}$ are \emph{internally disjoint} if 
$P$ and $Q$ do not intersect other than in these ordered ends.  {Note that the proof of this proposition uses Proposition~\ref{pathsincompletegraph}. As a consequence this proposition and each of the remaining results of this section, including the diameter lemma, require that $(k-\ell) \nmid k$.}

\begin{prop} \label{p(k,l)}
Suppose that $k \geq 3$ and that $k/2 \leq \ell \leq k-1$ is such that $(k-\ell) 
\nmid k$. Then there exists a $(4\ell-k+2)$-partite $k$-graph $\Pa(k, \ell)$ 
such that the following conditions hold.
\begin{itemize}
\item[(1)] $\Pa(k, \ell)$ is the union of $4\ell+1$ internally disjoint 
$\ell$-paths, each containing between $k^2\ell$ and $2k^5$ vertices, with 
identical ordered $\ell$-sets $T_1$ and $T_2$ as ordered ends (we refer to these 
as the ordered ends of $\Pa(k, \ell)$). In particular, $\Pa(k, \ell)$ contains at most 
$10k^6$ vertices.
\item[(2)] The vertex classes of $\Pa(k,\ell)$ are disjoint sets $V_w$, one for 
each vertex $w$ of $\W(k, \ell)$.
\item[(3)] Whenever $v_1 \in V_{w_1}, v_2 \in V_{w_2}, \dots, v_k \in V_{w_k}$ 
are such that $\{v_1, v_2, \dots, v_k\}$ is an edge of $\Pa(k, \ell)$, $\{w_1, 
\dots, w_k\}$ is an edge of $\W(k, \ell)$. Furthermore, let $X$ and $Y$ be the
ordered ends of $\W(k, \ell)$. 
Then the ordered ends of $\Pa(k, \ell)$ are contained in $\bigcup_{w \in X} 
V_{w}$ and $\bigcup_{w \in Y} V_{w}$ respectively.
\end{itemize}
\end{prop}  

\proof For every vertex $w$ of $\W(k, \ell)$, take a large vertex set $V_w$. 
Define $\W^*$ to have vertex set $V = \bigcup_{w \in \W(k, \ell)} V_w$, and 
edges precisely those $k$-sets of vertices which lie in sets corresponding to an 
edge of $\W(k, \ell)$. We shall construct $\Pa(k, \ell)$ to be a sub-$k$-graph 
of $\W^*$, with the ordered ends of $\Pa(k, \ell)$ in the sets $V_w$ 
corresponding to 
the ordered ends of $\W(k, \ell)$. Then $\Pa(k, \ell)$ will be a $(4\ell-k+2)$-partite 
$k$-graph which satisfies (2) and~(3). 

For each $1 \leq i \leq 2\ell - k +2$ let $e_i$ be the $i$th edge of $\W(k, 
\ell)$ as in the definition of $\W(k, \ell)$. Then for each $1 \leq i \leq 2\ell 
- k + 1 $ we know that $|e_i \cap e_{i+1}| = k-2$, and so we may choose $S_i$ to 
be an ordered set of $\ell$ vertices chosen from $\bigcup_{w \in e_i \cap 
e_{i+1}} V_w$. Also, let $S_0$ and $S_{2\ell-k+2}$ be ordered sets of $\ell$ 
vertices chosen from the $V_w$ corresponding to the ordered ends of $\W(k, \ell)$. So
$S_0$ and $S_{2\ell-k+2}$ are subsets of $\bigcup_{w \in e_1} V_w$, and
$\bigcup_{w \in e_{2\ell-k+2}} 
V_w$ respectively. We choose these sets $S_i$ to be disjoint and to contain at 
most one vertex from any one vertex class $V_w$. Then by 
Proposition~\ref{pathsincompletegraph}, for each $1 \leq i \leq 2\ell-k+2$ we 
can find an $\ell$-path from $S_{i-1}$ to $S_i$ in $\K[V_w: w \in e_i]$ which 
contains $k^2\ell(k-\ell)+k$ vertices. We do this so that the $\ell$-paths 
chosen only 
intersect in the appropriate $S_i$. Then the union of all of these $\ell$-paths 
is an $\ell$-path $P$ from $S_0$ to $S_{2\ell-k+2}$ with 
$$ k^2\ell \leq k^2\ell(k-\ell)+k \leq |P| \leq (2\ell-k+2)(k^2\ell(k-\ell)+k) 
\leq 2k^5.$$
In the same way we find another 
$4\ell$ $\ell$-paths from $S_0$ to $S_{2\ell-k+2}$, so that all $4\ell+1$ of the 
$\ell$-paths obtained are internally disjoint. Then the union of all of these 
$\ell$-paths is the $\Pa(k, \ell)$ we seek.
\endproof

Fix any such $\Pa(k, \ell)$, which we shall mean when we refer to $\Pa(k, \ell)$ 
in the rest of this paper. Also, let $S_1$ and $S_2$ be the ordered ends of 
$\Pa(k, \ell)$, so that $S_1$ and $S_2$ are disjoint ordered $\ell$-sets. Let 
$\Sa(k, \ell)$ be the complex with vertex set $S_1 \cup S_2$ and with edges 
being all subsets of $S_1$ and all subsets of $S_2$. Then since each of the 
$\ell$-paths which form $\Pa(k, \ell)$ contain at least $k^2\ell$ vertices, the 
complex $\Sa(k, \ell)$ is an induced subcomplex of the complex $\Pa(k, \ell)^\le$
corresponding to $\Pa(k, \ell)$, so under appropriate 
circumstances we will be able to use the extension lemma 
(Lemma~\ref{extensions_count}) to extend $\Sa(k, \ell)$ to $\Pa(k, \ell)$. This 
is the key to the following lemma, which states that for the values of $k$ and 
$\ell$ considered, almost all pairs of ordered $\ell$-sets of vertices of a 
sufficiently large $k$-graph of large minimum degree form the ordered ends of a 
copy of 
$\Pa(k, \ell)$.

\begin{lemma} \label{neardiameter}
Suppose that $k \geq 3$, that $k/2 \leq \ell \leq k-1$ is such that 
$(k-\ell)\nmid k$, and that $1/n \ll \beta \ll \mu, 1/k$. Let $\Hy$ be a 
$k$-graph of order $n$ with $\delta(\Hy) \geq \mu n$. Then there are at 
most~$\beta n^{2\ell}$ pairs of ordered $\ell$-sets~$S_1$ and~$S_2$ of vertices 
of $\Hy$ for which $\Hy$ does not contain a copy of $\Pa(k, \ell)$ with ordered 
ends $S_1$ and $S_2$. 
\end{lemma} 

\proof To prove this, we use hypergraph regularity. So introduce new constants
$$\frac{1}{n} \ll \frac{1}{r},\delta \ll c \ll \min\{\delta_k,1/t\} \leq 
\delta_k, \eta \ll d \ll \gamma \ll \beta 
\ll \mu.$$
We may assume that $t!$ divides $|\Hy|$, so apply the regularity lemma to $\Hy$, 
and let $V_1, \dots, V_{a_1}$ be the clusters of the partition obtained. As in 
Section~\ref{sec:rhgraph}, we say that an edge of $\Hy$ is useful if it lies in 
$\K_k(\hat{P}^{(k-1)})$ such that $\Hy$ is $(d', \delta_k, r)$-regular with 
respect to $\hat{P}^{(k-1)}$ for some $d' \geq d$. Let $\Hy'$ be the subgraph of 
$\Hy$ consisting of all useful edges. Note that no edge of $\Hy'$ 
contains 2 vertices from the same cluster. Then by Lemma~\ref{mostedgesgood}, 
at most $2d n^k$ edges of $\Hy$ are not useful, and so $d_{\Hy'}(S) \geq \mu n/2$ 
for all but at most~$\gamma n^{k-1}$ of the $(k-1)$-sets~$S$ of vertices 
of~$\Hy'$.

Let $C_1$ and $C_2$ be cells of the partition $\Part^{(\ell)}$ obtained from the 
regularity lemma. We say that $C_1$ and $C_2$ are \emph{connected} if $\Hy'$ 
contains a copy $\W$ of $\W(k, \ell)$ with ordered ends $A$ and $B$ such that $A \in 
C_1$, $B \in C_2$, and such that no two vertices of $\W$ lie in the same 
cluster. We shall first show that there are at most $\beta n^{2\ell}/2$ pairs 
$A$ and $B$ of ordered $\ell$-sets of vertices of $\Hy$ such that either 
\begin{itemize}
\item [(i)] at least one of $A$ and $B$ does not lie in a cell of 
$\Part^{(\ell)}$, or
\item [(ii)] the 
cells $C_A$ and $C_B$ of $\Pa^{(\ell)}$ which contain $A$ and $B$ respectively 
are not connected.
\end{itemize} 
Indeed, for (i) note that at most $\ell^2\frac{n}{a_1}n^{\ell-1} \leq 
\ell^2\eta n^\ell$ ordered $\ell$-sets of vertices of $\Hy$ do not
lie in $\cross_\ell$, and so there are at most $\ell^2\eta n^{2\ell}$
pairs $A$ and $B$ of ordered 
$\ell$-sets of vertices of $\Hy$ such that at least one of $A$ and $B$ does not 
lie in a cell of $\Part^{(\ell)}$. Similarly, for (ii) note that there are at most 
$\ell^2\eta n^{2\ell}$ pairs $A$ and $B$ of ordered $\ell$-sets such 
that the cells $C_A$ and $C_B$ of $\Pa^{(\ell)}$ which contain $A$ and $B$ 
respectively share at least one cluster. Finally, if the cells $C_A$ 
and $C_B$ of $\Pa^{(\ell)}$ which contain $A$ and $B$ respectively do not share 
any clusters, but are not connected, then $\Hy'$ must contain fewer 
than $\binom{4\ell-k+2}{2}\eta n^{2\ell-k+2}$ copies of $\W(k, \ell)$ with ordered
ends $A$ and $B$. So by 
Proposition~\ref{connectedclusters}, there are at most $\beta n^{2\ell}/3$ pairs 
$A$ and $B$ of ordered $\ell$-sets of vertices of $\Hy$ which lie in such pairs 
of cells of $\Part^{(\ell)}$. 

To prove the lemma, it is therefore sufficient to show that  there are at most 
$\beta n^{2\ell}/2$ pairs $S_1, S_2$ of ordered $\ell$-sets of vertices of $\Hy$ 
such that $C_{S_1}$ and $C_{S_2}$ are connected cells of $\Part^{(\ell)}$ but 
$S_1$ and $S_2$ do not form the ordered ends of a copy of $\Pa(k,\ell)$ in 
$\Hy$. So suppose cells $C_1$ and $C_2$ of $\Pa(k,\ell)$ are connected. Then 
there is a copy $\W$ of $\W(k,\ell)$ in $\Hy'$ with ordered ends $A$ and $B$ such that 
$A \in C_1$, $B \in C_2$, and such that no two vertices of $\W$ lie in the same 
cluster. Since every edge of $\Hy'$ is a useful edge, for each edge $e \in 
E(\W)$ the polyad $\hat{P}^{(k-1)}(e)$ of $e$ is such that $\Hy$ is $(d', \delta_k, r)$-regular 
with respect to $\hat{P}^{(k-1)}(e)$ for some $d' \geq d$. Then these 
polyads `fit together'. By this we mean that if edges $e$ and $e'$ of $\W$ 
intersect in $q$ vertices, then
$$
\left(\bigcup_{i=1}^{k-1} \hat{P}^{(i)}(e)\right)\cap
\left(\bigcup_{i=1}^{k-1} \hat{P}^{(i)}(e')\right)=\bigcup_{i=1}^{q} \hat{P}^{(i)}(e\cap e'),
$$
i.e.~the intersection of the $(k-1,k)$-complexes corresponding to~$e$ and~$e'$ is
the $(q,q)$-complex corresponding to $e\cap e'$.
Therefore we can define $\Hy^*$ to be the $(k, 4\ell-k+2)$-complex obtained from 
the $(k-1, 4\ell-k+2)$ complex $\bigcup_{e \in E(\W)} \bigcup_{i=1}^{k-1} 
\hat{P}^{(i)}(e)$ by adding $E(\Hy) \cap \bigcup_{e \in E(\W)} 
\K(\hat{P}^{(k-1)}(e))$ as the `$k$th level'. Then $\Hy^*$ is a  
$(\mathbf{d},\delta_k,\delta,r)$-regular $(k, 4\ell - k +2)$-complex, where 
$\mathbf{d} = (d', d_{k-1}, \dots, d_2)$ and $(d_{k-1}, \dots, d_2)$ is as in 
the definition of a $(\eta, \delta, t)$-equitable family of partitions.
(Here we may assume a common density~$d'$ for the $k$th level by applying
the slicing lemma (\cite{RS}, Proposition~22) if necessary.)
Furthermore, by construction $\Hy^*$ respects the partition of the complex
$\W(k, \ell)^\le$ corresponding to $\W(k, \ell)$, and so property~(3) of
Proposition~\ref{p(k,l)} implies
that $\Hy^*$ also respects the partition of $\Pa(k, \ell)^\le$.
Let $S_1$ and $S_2$ be disjoint ordered $\ell$-sets lying in the cells $C_1$ 
and $C_2$ of  $\Part^{(\ell)}$ respectively. Then $S_1 \cup S_2$ is the vertex 
set of a labelled copy $\Sa$ of $\Sa(k, \ell)$ in $\Hy^*$. So by Lemma~\ref{extensions_count}, for all but at most $\beta |C_1||C_2|/2$ choices of $S_1 
\in C_1$ and $S_2 \in C_2$ we can extend the labelled complex $\Sa$ to at least 
one labelled partition respecting copy of 
$\Pa(k, \ell)$ with ordered ends $S_1$ and $S_2$. Summing over all $C_1$ and 
$C_2$, we find that there are at most $\beta n^{2\ell}/2$ ordered $\ell$-sets 
$S_1$ and $S_2$ of vertices of $\Hy$ which lie in connected cells of 
$\Pa^{(\ell)}$ and which cannot be extended to a labelled partition respecting 
copy of $\Pa(k, \ell)$, completing the proof.
\endproof

We can now prove the following corollary, the diameter lemma we were aiming for. 
The idea behind this is that if $S$ and $T$ are ordered $\ell$-sets in a large 
$k$-graph $\Hy$ of large minimum degree, then there are many ordered $\ell$-sets 
$S'$ and $T'$ such that $\Hy$ contains $\ell$-paths from $S$ to $S'$ and $T$ to 
$T'$. So by the previous lemma, at least one such pair $S'$ and $T'$ will form 
the ordered ends of a copy of $\Pa(k, \ell)$, and then combining these 
$\ell$-paths we 
will obtain an $\ell$-path from $S$ to $T$.

\begin{coro}[Diameter lemma] \label{diameter}
Suppose that~$k \geq 3$, that~$1 \leq \ell \leq k-1$ is such that that~$(k-\ell) 
\nmid k$, and that~$1/n \ll \mu, 1/k$. Let~$\Hy$ be a~$k$-graph of order~$n$ 
with~$\delta(\Hy) \geq \mu n$. Then for any two disjoint ordered~$\ell$-sets $S$ 
and $T$ of vertices of~$\Hy$, there exists an~$\ell$-path~$P$ in~$\Hy$ from~$S$ 
to~$T$ such that~$P$ contains at most~$8k^5$ vertices.
\end{coro}

\proof Recall that if~$\ell < k/2$ we can find such an $\ell$-path consisting of 
just one single edge, so we may assume that~$\ell \geq k/2$. Introduce constants 
$\beta,\beta'$ such that $1/n \ll \beta' \ll \beta \ll \mu, 1/k$. Let $A$ be an 
arbitrary ordered $\ell$-set of vertices of $\Hy$, and let $X$ be an arbitrary 
set of 
$3\ell$ vertices which is disjoint from $A$. We begin by showing that there are 
many ordered $\ell$-sets $B$ such that $\Hy$ contains an $\ell$-path $P$ from 
$A$ to $B$ having at most $3\ell$ vertices, none of which are from $X$. To 
show this, we shall demonstrate how a vertex sequence of $P$ may be chosen, and 
then count the number of choices. 

Since $A$ will be an ordered end of $P$, we begin the vertex sequence of $P$ 
with the ordered $\ell$-set $A$. We then arbitrarily choose any $k-\ell-1$ 
vertices of $\Hy$ to add to the sequence. To finish the sequence, we repeatedly 
make use of the fact that $\delta(\Hy) \geq \mu n$. More precisely, we repeat 
the following step: let $V$ be the set of the final $k-1$ vertices of the 
current vertex sequence. Then $d_\Hy(V) \geq \mu n$, and so there are at least 
$\mu n - 6\ell $ vertices which together with $V$ form an edge of $\Hy$ and which are 
not in the vertex sequence constructed thus far or in $X$. Choose $v$ to be one 
of these vertices, and append it to the vertex sequence. We stop as soon as the number 
$r$ of vertices in the sequence satisfies $r > 2 \ell$ and $r \equiv k$ 
(modulo $(k-\ell)$), so in particular $r \leq 3 \ell$. Let $B$ be the ordered 
set consisting of the last $\ell$ vertices of the sequence. Then $\Hy$ contains 
an $\ell$-path $P$ with this vertex sequence, and $P$ is therefore an 
$\ell$-path of order at most $3\ell$ from $A$ to $B$ which does not contain any 
vertex of $X$. There are at least $(\mu n -6\ell)^{r-\ell}$ vertex sequences we 
could have chosen, and hence there are at least $(\mu n 
-6\ell)^{r-\ell}/n^{r-2\ell} > \beta n^\ell$ possibilities for an ordered 
$\ell$-set $B$ such that there is an $\ell$-path from~$A$ to~$B$ in~$\Hy$, not
containing any vertex of $X$.

Now, let $S$ and $T$ be the two ordered $\ell$-sets of vertices of $\Hy$ given 
in the statement of the corollary. Then there are at least $\beta n^{\ell}$ 
ordered $\ell$-sets $S'$ of vertices of $\Hy$ such that there exists an 
$\ell$-path $P_1$ from $S$ to $S'$ in $\Hy$, which contains at most $3\ell$ 
vertices and such that $V(P_1) \cap T = \emptyset$. Likewise for each such 
choice of $S'$ and $P_1$, there are at least $\beta n^{\ell}$ ordered 
$\ell$-sets $T'$ of vertices of $\Hy$ such that there exists an $\ell$-path 
$P_2$ from $T$ to $T'$ of order at most $3 \ell$ in $\Hy$ and such that $V(P_2) 
\cap V(P_1) = \emptyset$. By Lemma~\ref{neardiameter}, at most $\beta' n^{2 
\ell}$ of these pairs $S', T'$ do not form ordered ends of a copy of $\Pa(k, 
\ell)$ in $\Hy$. Since $\beta' \ll \beta$ we may therefore choose such a pair $S', T'$ 
such that $S'$ and $T'$ are ordered ends of a copy of $\Pa(k, \ell)$ in $\Hy$. 
Then there are at least $4 \ell+1$ internally disjoint $\ell$-paths of order at 
most $2k^5$ from $S'$ to $T'$ in $\Hy$. At most $4 \ell$ of these 
$\ell$-paths 
contain any vertex from $V(P_1)\sm S'$ or $V(P_2)\sm T'$, and so we may choose 
an $\ell$-path $Q$ from $S'$ to $T'$ in $\Pa(k, \ell) \subseteq \Hy$ of order 
at most $2k^5$ which contains no vertex from $V(P_1)\sm S'$ or $V(P_2)\sm 
T'$. Then $P_1QP_2$ is the $\ell$-path from $S$ to $T$ of order at most 
$2k^5+6\ell \leq 8k^5$ we seek. \endproof

\section{Absorbing Path Lemma} \label{sec:absorb}

Let~$\Hy$ be a $k$-graph, and let~$S$ be a set of~$k-\ell$ vertices of~$\Hy$. Recall
that an $\ell$-path~$P$ in~$\Hy$ with ordered ends~$P^{beg}$ and~$P^{end}$ 
is \emph{absorbing} for~$S$ if~$P$ does not contain any vertex 
of~$S$, and~$\Hy$ contains an $\ell$-path~$Q$ with the same ordered 
ends~$P^{beg}$ and~$P^{end}$, where $V(Q) = V(P) \cup S$. This means that if~$P$ 
is a section of an $\ell$-path~$P^*$ which does not contain any vertices of~$S$, 
then we can `absorb' the vertices of~$S$ into~$P^*$ by replacing~$P$ with~$Q$. 
$P^*$ will still be an~$\ell$-path after this change as~$P$ and~$Q$ have the 
same ordered ends. Similarly, we say that an $\ell$-path $P$ in $\Hy$ with 
ordered ends~$P^{beg}$ and~$P^{end}$ can \emph{absorb} a collection $S_1, \dots, 
S_r$ of $(k-\ell)$-sets of vertices of $\Hy$ if~$P$ does not contain any vertex 
of~$\bigcup_{i=1}^r S_i$, and~$\Hy$ contains an $\ell$-path~$Q$ with the same 
ordered ends~$P^{beg}$ and~$P^{end}$, where $V(Q) = V(P) \cup \bigcup_{i=1}^r 
S_i$. The reason we absorb $(k-\ell)$-sets is that the order of 
an $\ell$-path must be congruent to $k$, modulo $k-\ell$. The next result 
describes the absorbing path as a $k$-graph, which we shall use to 
absorb a set~$S$.  {Note that the proof of this proposition uses Proposition~\ref{pathsincompletegraph}. As a consequence, this proposition and each of the remaining results of this section, including the absorbing path lemma, require that $(k-\ell) \nmid k$.}

\begin{prop}\label{absorbingpathsegment}
Suppose that $k \geq 3$, and that $1 \leq \ell \leq k-1$ is such that 
$(k-\ell)\nmid k$. Then there is a $k$-partite $k$-graph $\AP(k,\ell)$ with the 
following properties.
\begin{itemize}
\item[(1)] $|\AP(k, \ell)| \leq k^4.$
\item[(2)] The vertex set of $\AP(k, \ell)$ consists of two disjoint sets $S$ 
and $X$ with $|S| = k-\ell$.
\item[(3)] $\AP(k, \ell)$ contains an $\ell$-path $P$ with vertex set $X$ and 
ordered ends $P^{beg}$ and $P^{end}$.
\item[(4)] $\AP(k, \ell)$ contains an $\ell$-path $Q$ with vertex set $S \cup X$ 
and ordered ends $P^{beg}$ and $P^{end}$.
\item[(5)] No edge of $\AP(k, \ell)$ contains more than one vertex of $S$.
\item[(6)] No vertex class of $\AP(k, \ell)$ contains more than one vertex of 
$S$.
\end{itemize}
\end{prop} 

\proof Let $V_1, \dots, V_k$ be disjoint vertex sets of size $k\ell(k-\ell)+1$. Let $S$ be an
ordered $(k-\ell)$-set such that for each $1 \leq i \leq k-\ell$, the $i$th vertex of $S$
lies in $V_{\ell+i}$. Let $P$ be an $\ell$-path in $\K[V_1, \dots, V_k]$ with ordered ends
$P^{beg}$ and $P^{end}$ such that both $P^{beg}$ and $P^{end}$ contain at most one
vertex from each~$V_i$ and such that $V(P)=(V_1\cup\dots\cup V_k)\setminus S$.
(One can easily choose such a~$P$ if for all $j=1,\dots,|P|$ one chooses the $j$th
vertex of~$P$ in the $V_i$ for which $j\equiv i$ modulo~$k$.)
Then $V(P)\cup S=V_1\cup\dots\cup V_k$.
Thus we can apply Proposition~\ref{pathsincompletegraph} to obtain an $\ell$-path
$Q$ from $P^{beg}$ to $P^{end}$ in $\K[V_1, \dots, V_k]$ such that $V(Q)=V(P)\cup S$.
By swapping some vertices in~$S$ with some vertices in $V(Q)\setminus S$ (lying in the same~$V_i$)
if necessary we can ensure that the vertices in~$S$ are distributed in
such a way that in some vertex sequence of~$Q$ they have distance at least~$k$ from each other.
(This ensures~(5).) We can now take $\AP(k, \ell):=P\cup Q$.\endproof

Fix an~$\AP(k, \ell)$ satisfying Proposition~\ref{absorbingpathsegment}, which 
we shall refer to simply as~$\AP(k, \ell)$ for the rest of this section. Let 
$b(k, \ell) :=|\AP(k, \ell)|-k+\ell$, so that~$b(k, \ell)$ is the number of 
vertices of the $\ell$-path~$P$ in the definition of~$\AP(k ,\ell)$. 

Now, given a $(k-\ell)$-set~$S$ of vertices of~$\Hy$, we can think of~$S$ as a 
labelled $(k, k)$-complex with no $i$-edges for any~$i \geq 2$. We will apply 
the extension lemma (Lemma~\ref{extensions_count}) to deduce that for most such 
$(k-\ell)$-sets $S$, there are many labelled copies of~$\AP(k, \ell)^\leq$ 
extending~$S$ in $\Hy$, which will imply that~$\Hy$ contains many absorbing paths 
for these sets~$S$.

Suppose that~$\Hy$ is a $k$-graph on~$n$ vertices, and that~$c$ is a positive 
constant. We say that a $(k-\ell)$-set~$S$ of vertices of~$\Hy$ is 
\emph{$c$-good} if~$\Hy$ contains at least~$c n^{b(k, \ell)}$ absorbing paths
for~$S$, each on $b(k, \ell)$ vertices. $S$ is \emph{$c$-bad} if it is not $c$-good. 
The next lemma states that for the 
values of~$k$ and~$\ell$ we are interested in, and any small~$c$, if~$\Hy$ is 
sufficiently large and has large minimum degree, then almost all 
$(k-\ell)$-sets~$S$ of vertices of~$\Hy$ are $c$-good.

\begin{lemma} \label{fewbadtuples}
Suppose that $k \geq 3$, that $1 \leq \ell \leq k-1$ is such that $(k-\ell)\nmid 
k$, and that $1/n \ll c \ll \gamma\ll \mu, 1/k$. Let $\Hy$ be a $k$-graph on $n$ 
vertices 
such that $\delta(\Hy) \geq \mu n$. Then at most $\gamma n^{k-\ell}$ sets $S$ of 
$k-\ell$ vertices of $\Hy$ are $c$-bad.
\end{lemma}

\proof Let~$b=b(k, \ell)$, and introduce new constants 
$$ \frac{1}{n} \ll \frac{1}{r},\delta \ll c \ll
\min\{\delta_k,1/t\} \leq \delta_k, \eta \ll d \ll \gamma. $$
We may assume that~$t!$ divides~$|\Hy|$, so apply the regularity lemma to~$\Hy$, 
and let $V_1, \dots, V_{a_1}$ be the clusters of the partition obtained. Let $m 
= n/a_1$ be the size of each of these clusters. Form the reduced 
$k$-graph~$\R$ on these clusters as defined in Section~\ref{sec:rhgraph}. 

We begin by showing that almost all sets of~$k-\ell$ vertices of~$\Hy$ are 
contained in clusters which lie in a common edge of $\R$. More precisely, for 
all but at most $\gamma n^{k-\ell}/2$ sets $\{v_1, \dots, v_{k-\ell}\}$ 
of~$k-\ell$ vertices of~$\Hy$ we can choose clusters $V_{i_1}, \dots, V_{i_k}$ 
such that $v_j \in V_{i_j}$ for each $1 \leq j \leq k-\ell$ and such that 
$\{V_{i_1}, \dots, V_{i_k}\}$ forms an edge of~$\R$. 
Indeed, by Lemma~\ref{minreduced}, $d_\R(S) \geq 1$ for all but at most $\gamma 
a_1^{k-\ell}/3$ `neighbourless' sets~$S$ of~$k-\ell$ clusters. At most $\eta 
n^{k-\ell} \ll \gamma n^{k-\ell}$ sets~$T$ of~$k-\ell$ vertices of~$\Hy$ do not 
lie in $\cross_{k-\ell}$. But if $T \in \cross_{k-\ell}$, then unless the set $S$ of 
clusters containing the vertices of $T$ is one of the at most $\gamma 
a_1^{k-\ell}/3$ `neighbourless' sets of~$k-\ell$ clusters (which will be the case 
for at most $\gamma n^{k-\ell}/3$ sets of~$k-\ell$ vertices of~$\Hy$), there is 
an edge of~$\R$ containing~$S$ as required.

Now, suppose that $V_{i_1}, \dots, V_{i_k}$ are clusters which form an edge of~$\R$.
Note that there are $m^{k-\ell}$ sets $\{v_1, \dots, 
v_{k-\ell}\}$ such that $v_j \in V_{i_j}$ for each $1 \leq j \leq k-\ell$. Since 
$e=\{V_{i_1}, \dots, V_{i_k}\}$ is an edge of~$\R$, we may define the 
complex $\Hy^*$ corresponding to~$e$ as in the paragraph after the 
statement of the extension lemma (Lemma~\ref{extensions_count}).
Then $\Hy^*$ satisfies the conditions of the 
extension lemma  (with $\gamma/2$ and $k$ playing the roles
of $\beta$ and $s$), and respects the partition of $\AP(k, \ell)$. Let $S$ be an 
ordered set of size $k-\ell$, which we may view as a labelled $(k, k)$-complex 
with no $j$-edges for $j \geq 2$. Then by Lemma~\ref{extensions_count} applied with $cb(k, \ell)! a_1^{b(k, \ell)}$ in place of $c$, all but 
at most $\gamma m^{k-\ell}/2$ ordered sets $S' = \{v_1, \dots, v_{k-\ell}\}$ such that 
$v_j \in V_{i_j}$ for each $j$ (these are the labelled copies of $S$) are 
extendible to at least $cb(k, \ell)!a_1^{b(k, \ell)}m^{b(k,\ell)}$ labelled partition-respecting
copies of $\AP(k,\ell)$ in $\Hy$. This is where we use property (5) of 
Proposition~\ref{absorbingpathsegment} \textendash~it ensures that the complex 
$S$ is an induced subcomplex of $\AP(k, \ell)$. For each copy $C$ of 
$\AP(k, \ell)$, $C-S'$ is an absorbing path for $S'$ on $b(k,\ell)$ vertices, and
so $\Hy^*$ (and therefore $\Hy$) contains at least 
$cn^{b(k, \ell)}$ absorbing paths on $b(k, \ell)$ vertices for~$S'$. So at most
$\gamma m^{k-\ell}/2$ such sets~$S'$ are $c$-bad.
 
Recall that the number of $(k-\ell)$-sets of vertices 
of~$\Hy$ which do not lie in distinct clusters corresponding to an edge of~$\R$
is at most $\gamma n^{k-\ell}/2$. Summing over all sets of $k-\ell$ clusters, we 
see that at most $\gamma n^{k-\ell}/2$ of the $(k-\ell)$-sets which do lie in 
distinct clusters corresponding to an edge of~$\R$ are $c$-bad. Thus at most 
$\gamma n^{k-\ell}$ sets of $k-\ell$ vertices of~$\Hy$ are $c$-bad, completing 
the proof.\endproof

We are now in a position to prove the main lemma of this section. It states 
that for any positive $c$, if $\Hy$ is a sufficiently large $k$-graph of large 
minimum degree, then we can 
find an $\ell$-path in~$\Hy$ which contains a small proportion of the vertices 
of~$\Hy$, includes all vertices of~$\Hy$ which lie in many $c$-bad 
$(k-\ell)$-sets 
and can absorb any small collection of $c$-good $(k-\ell)$-sets of vertices 
of~$\Hy$.

\begin{lemma}[Absorbing path lemma] \label{absorb}
Suppose that $k \geq 3$, that $1 \leq \ell \leq k-1$ is such that $(k-\ell)\nmid k$, 
and that $1/n \ll \alpha \ll c \ll \gamma\ll \mu, 1/k$. Let $\Hy$ be a $k$-graph 
of order $n$ with $\delta(\Hy) \geq \mu n$. Then $\Hy$ contains an $\ell$-path $P$
on at most $\mu n$ vertices such that the following properties hold:
\begin{itemize}
\item[(1)] Every vertex of $\Hy - V(P)$ lies in at most $\gamma n^{k-\ell-1}$ 
$c$-bad $(k-\ell)$-sets. 
\item[(2)] $P$ can absorb any collection of at most $\alpha n$ disjoint $c$-good 
$(k-\ell)$-sets of vertices of~$\Hy - V(P)$.
\end{itemize}
\end{lemma}

\proof Let $b := b(k, \ell)$, and choose a family $\T$ of ordered $b$-sets of vertices 
of $\Hy$ at random by including each ordered $b$-set $T$ into $\T$ with probability
$c^2n^{1-b}$, independently of all other ordered $b$-sets. 
Now, for any $c$-good set $S$ of $k-\ell$ vertices of~$\Hy$, the expected number of
$T\in \T$ for which~$\Hy$ contains an absorbing path for $S$ with $T$ as a vertex
sequence is at least~$c^3n$, by the definition of a $c$-good set. So by a standard 
Chernoff bound, with probability $1-o(1)$, for every $c$-good $(k-\ell)$-set 
$S$ of vertices of $\Hy$ the number of such ordered $b$-sets $T\in \T$ is
at least $c^3n/2$. Furthermore, with probability $1-o(1)$ we 
have $|\T| \leq 2c^2n$. 
The expected number of ordered pairs $T, T'$ in $\T$ which intersect (i.e.~for
which the corresponding unordered sets intersect) is at most 
$(c^2n^{1-b})^2b^2n^{2b-1}= c^4b^2n$. So with probability at least
1/2 the number of such pairs is at most $2c^4 b^2n$. Thus we may fix an outcome of our 
random selection of $\T$ such that all of these events hold. 

Delete from $\T$ every $T \in \T$ which intersects any other $T' \in \T$. Also 
delete from $\T$ every $T \in \T$ which is not a vertex sequence of an 
absorbing path for some $c$-good 
$(k-\ell)$-set $S$ of vertices of $\Hy$. Let $T_1, \dots, T_q$ be the remaining 
members of $\T$. So $q \leq 2c^2 n$, and for each $1 \leq i \leq q$ we can 
choose an $\ell$-path~$P_i$ in $\Hy$ with vertex sequence $T_i$ which is absorbing for some 
such $S$. Then all the $\ell$-paths $P_i$ are 
disjoint, and for every $c$-good $(k-\ell)$-set $S$ of vertices of $\Hy$ at least 
$c^3n/2 - 2c^4b^2n \geq \alpha n$ of the $\ell$-paths $P_i$ are absorbing.

Let $X$ be the set of vertices of $\Hy$ which are not contained in any $P_i$ and 
which lie in more than~$\gamma n^{k-\ell-1}$ $c$-bad~$(k-\ell)$-sets. Then $|X| 
\leq \gamma n$, since by Lemma~\ref{fewbadtuples} there are at most $\gamma^2 
n^{k-\ell}/(k-\ell)$ $c$-bad 
$(k-\ell)$-sets in total. We shall use the minimum degree condition on~$\Hy$ to 
greedily construct an~$\ell$-path $P_0$ containing all vertices in $X$ and not 
intersecting the previous paths~$P_i, 1 \leq i \leq q$. Then if we incorporate 
each of the $P_i$ ($0 \leq i \leq q$) into the $\ell$-path $P$ we are constructing, 
conditions~(1) and~(2) of the lemma will be satisfied. So let $X'$ be a set of $k-1$ vertices of 
$X$. Then $d_\Hy(X') \geq \mu n$ by the minimum degree condition on $\Hy$. Since 
$|\bigcup_{i=1}^{q} P_i| < \mu n$, we may choose a vertex $x \in V(\Hy) \sm 
\bigcup_{i=1}^{q} V(P_i)$ which together with $X'$ forms an edge of $\Hy$. Then $X' 
\cup \{x\}$ is the first edge of $P_0$. We then greedily extend $P_0$ as 
follows. Let $X''$ be the set of the final $\ell$ vertices of the vertex 
sequence of $P_0$. Add to $X''$ any $k-1-\ell\ge 1$ vertices from $X$ not yet 
contained in $P_0$. Then $d_\Hy(X'') \geq \mu n$, and so we may choose a vertex 
$y$ of $\Hy$ which is not in $\bigcup_{i=1}^{q} P_i$ nor already contained in 
$P_0$. We then extend $P_0$ by the edge $X'' \cup \{y\}$. At the end of this 
process we obtain an $\ell$-path $P_0$ which is disjoint from all the $P_i$ ($i=1,\dots q$), 
which contains every vertex of $X$, and which satisfies $|V(P_0)| \leq 2\gamma n$.
Let $P_i^{beg}$ and $P_i^{end}$ be ordered ends of $P_i$ for each $0 \leq i \leq q$.

To complete the proof, we now use the diameter lemma (Corollary~\ref{diameter}) to 
greedily join each ordered $\ell$-set $P_i^{end}$ to the ordered $\ell$-set
$P_{i+1}^{beg}$ by an 
$\ell$-path $P'_i$, such that $P'_i$ intersects $P_i$ and $P_{i+1}$ only in the 
sets $P_{i+1}^{beg}$ and $P_i^{end}$ and does not intersect any other $P_j$ or 
any previously chosen $P'_j$. More precisely, suppose we have chosen such $P'_0, 
\dots, P'_{i-1}$. Let $\Hy'$ be the $k$-graph obtained from~$\Hy$ by removing
all the vertices in $P_0,\dots,P_q$ and all the vertices in $P'_0,\dots,P'_{i-1}$
and then adding back $P_i^{end}$ and $P_{i+1}^{beg}$.
Then $\delta(\Hy') \geq \mu n/2$, and so we may apply Corollary~\ref{diameter} to 
find an $\ell$-path $P'_i$ in $\Hy'$ from $P_i^{end}$ to $P_{i+1}^{beg}$ 
containing at most $8k^5$ vertices.
Having found these $\ell$-paths, the absorbing path~$P^*$ is the $\ell$-path 
$P_0 P'_0 P_1 P'_1 P_2 \dots P_{q-1} P'_{q-1} P_q$.
\endproof

\section {Path Cover Lemma}\label{sec:pathcover}

\begin{lemma}[Path cover lemma]\label{pathcover}
Suppose $k \geq 3$, that $1 \leq \ell \leq k-1$, and that $1/n \ll 1/D \ll \eps \ll \mu,1/k$. 
Let $\Hy$ be a $k$-graph of order $n$ with $\delta(\Hy) \geq (\frac{1}{\lceil 
\frac{k}{k-\ell} \rceil (k-\ell)}+\mu)n$. Then $\Hy$ contains a set of at most 
$D$ disjoint $\ell$-paths covering all but at most $\eps n$ vertices of~$\Hy$.
\end{lemma}

Note that the condition $(k-\ell) \nmid k$ is not needed for this lemma. Let 
\begin{equation}\label{defa}
a:=\left\lceil\frac{k}{k-\ell}\right\rceil (k-\ell)
\end{equation}
and let $\F_{k,\ell}$ be the $k$-graph whose vertex set is the disjoint union of 
sets $A_1,\dots,A_{a-1}$ and $B$ of size $k-1$ and whose edges are all the 
$k$-sets of the form $A_i\cup \{b\}$ (for all $i=1,\dots,a-1$ and all $b\in 
B$). An \emph{$\F_{k,\ell}$-packing} in a $k$-graph $\R$ is a collection of pairwise 
vertex-disjoint copies of $\F_{k,\ell}$ in~$\R$.

The idea of the proof of the path cover lemma is to apply the regularity lemma 
to~$\Hy$ in order to obtain a reduced $k$-graph~$\R$. Recall that by 
Lemma~\ref{minreduced} the minimum degree of $\Hy$ is almost inherited by~$\R$. 
So we can use the following lemma (Lemma~\ref{Fklpacking}) to obtain an almost 
perfect $\F_{k,\ell}$-packing in~$\R$. Consider any copy~$\F$ of $\F_{k,\ell}$ 
in this packing. We will repeatedly apply the embedding lemma 
(Lemma~\ref{emblemma}) to the 
sub-$k$-graph $\Hy(\F)$ of~$\Hy$ corresponding to~$\F$ to obtain a bounded number of 
$\ell$-paths which cover almost all vertices of~$\Hy(\F)$. Doing this for all 
the copies of $\F_{k,\ell}$ in the $\F_{k,\ell}$-packing of~$\R$ will give a set 
of $\ell$-paths as required in Lemma~\ref{pathcover}.

\begin{lemma}\label{Fklpacking}
Suppose that $k \geq 3$, that $1 \leq \ell \leq k-1$, and that $1/n \ll \theta 
\ll \varepsilon \ll 1/k$. Let $\Hy$ be a $k$-graph of order $n$ such that $d(S) \geq 
n/a$ for all but at most $\theta n^{k-1}$ sets $S\in \binom{V(\Hy)}{k-1}$, where 
$a$ is as defined in~(\ref{defa}). Then $\Hy$ contains a $\F_{k,\ell}$-packing 
covering all but at most $(1-\varepsilon)n$ vertices.
\end{lemma}

We omit the proof of this lemma. It was first proved in~\cite{KO} for the case 
when $k=3$ and $\ell=1$. A proof for the case when $\ell<k/2$ can be found 
in~\cite{HS}. The general case can be proved similarly, see~\cite{Richard} for
details.

\begin{lemma}\label{pathcolouring}
Let $P$ be an $\ell$-path on $n$ vertices and let $a$ be as defined 
in~(\ref{defa}). Then there is a $k$-colouring of $P$ with colours $1,\dots,k$ 
such that colour $k$ is used $n/a\pm 1$ times and the sizes of all other 
colour classes are as equal as possible.
\end{lemma}

\proof
Let $x_1,\dots,x_n$ be a vertex sequence of $P$. Colour vertices 
$x_k,x_{k+a},x_{k+2a},\dots$ with colour $k$ and remove these vertices from the 
sequence $x_1,\dots,x_n$. Colour the remaining vertices in turn with colours 
$1,\dots,k-1$ as follows. Colour the first vertex with colour~$1$. Suppose that 
we just coloured the $i$th vertex with some colour~$j$. Then we colour the next 
vertex with colour $j+1$ if $j\le k-2$ and with colour~$1$ if $j=k-1$. To show 
that this yields a proper colouring, it suffices to show that every edge of~$P$ 
contains some vertex of colour~$k$. Clearly this holds for the first edge $e_1$ 
of $P$ and for all edges intersecting $e_1$ (since $x_k$ lies in all those 
edges). 
Note that the first vertex of the $i$th edge $e_i$ of $P$ is $x_{f(i)}$, where
$f(i)=(i-1)(k-\ell)+1$. Also note that $i^*:=\lceil \frac{k}{k-\ell} \rceil+1$
is the smallest integer so that $f(i^*)>k$.
In other words, the $i^*$th edge $e_{i^*}$ of $P$ is the first edge which does not contain
$x_k$. But the vertices of $e_{i^*}$ are $x_{a+1},\dots,x_{a+k}$.
So $e_{i^*}$ as well as all succeeding edges which 
intersect $e_{i^*}$ contain a vertex of colour~$k$ (namely~$x_{a+k}$). Continuing in 
this way gives the claim.
\endproof

\removelastskip\penalty55\medskip\noindent{\bf Proof of Lemma~\ref{pathcover}. }
Choose new constants such that
$$\frac{1}{n} \ll \frac{1}{D}\ll \frac{1}{r},\delta, c \ll \min\{\delta_k,1/t\} \leq \delta_k, 
\eta \ll d\ll \theta \ll \varepsilon.$$
We may assume that $t!|n$, so apply Theorem~\ref{reglemma} (the regularity lemma) 
to $\Hy$, and let $V_1, \dots, 
V_{a_1}$ be the clusters of the partition obtained. Let $m = n/a_1$ be the 
size of each of these clusters. Form the reduced $k$-graph $\R$ on these 
clusters as discussed in Section~\ref{sec:rhgraph}. Lemmas~\ref{minreduced} 
and~\ref{Fklpacking} together imply that~$\R$ has a 
$\F_{k,\ell}$-packing~$\mathcal{A}$ covering all but at most 
$\eps n/2$ vertices of~$\R$. Consider any copy $\F$ of $\F_{k,\ell}$ 
in this packing. Our aim is to cover almost all vertices in the clusters 
belonging to~$\F$ by a bounded number of disjoint $\ell$-paths.

So let $A_1,\dots,A_{a-1}$ and $B$ be $(k-1)$-element subsets of $V(\F)$ as in 
the definition of $\F_{k,\ell}$. So the edges of $\F_{k,\ell}$ are all the 
$k$-tuples of the form $A_i\cup \{b\}$ for all $i=1,\dots,a-1$ and all $b\in B$. 
Pick $b\in B$ and consider the edge $A_1\cup \{b\}=:e$. Let $\V$ be the set of 
all clusters corresponding to vertices in~$A_1$ and let $V_b$ be the cluster 
corresponding to $b$. Define the complex $\Hy^*$ corresponding to the edge~$e$ 
as in the paragraph after the statement of the extension lemma 
(Lemma~\ref{extensions_count}). Then Lemma~\ref{pathcolouring} and the embedding 
lemma (Lemma~\ref{emblemma} applied to $\Hy^*$) together imply that the 
sub-$k$-graph of $\Hy$ spanned by the vertices in~$V_b\cup \bigcup_{V\in \V} V$ 
contains an $\ell$-path $P_1$ on $acm/(a-1)$ vertices which intersects each cluster
from $\V$ in $cm/(k-1)\pm 1$ vertices and $V_b$ in $cm/(a-1) \pm 1$ vertices.
Lemma~\ref{restriction} implies that the subcomplex of $\Hy^*$ obtained by deleting 
the vertices of~$P_1$ is still $({\bf d},\sqrt{\delta_k},\sqrt{\delta},r)$-regular.
So we can find another $\ell$-path 
$P_2$ which is disjoint from $P_1$ and intersects each cluster from $\V$ in 
$cm/(k-1)\pm 1$ vertices and $V_b$ in $cm/(a-1) \pm 1$ vertices. We do this until we 
have used about $m/(k-1)$ vertices in each cluster from~$\V$. So we have found 
$1/c$ disjoint $\ell$-paths. Now we pick $b'\in B\setminus\{b\}$ and argue as 
before to get $1/c$ disjoint $\ell$-paths, such that each of them intersects 
(the remainder of) each cluster from $\V$ in $cm/(k-1)\pm 1$ vertices and 
$V_{b'}$ in $cm/(a-1) \pm 1$ vertices. We do this for all the $k-1$ vertices 
in~$B$. However, when considering the last vertex $b''$ of~$B$, we stop as soon as one 
of the subclusters from~$\V$ has size less than $\varepsilon m/4a$ (and thus all the 
other subclusters from~$\V$ have size at most $\varepsilon m/2a$) since we need 
to ensure that the subcomplex of $\Hy^*$ restricted to the remaining subclusters 
is still $({\bf d},\sqrt{\delta_k},\sqrt{\delta},r)$-regular. So in total we 
have chosen close to $(k-1)/c$ disjoint $\ell$-paths covering all but at most 
$\varepsilon m/2a$ vertices in each cluster from~$\V$ and covering between 
$m/(a-1) - \varepsilon m/2a$ and $m/(a-1)$ vertices in each 
cluster $V_b$ with $b\in B$. We now repeat this process for 
each of $A_2,\dots,A_{a-1}$ in turn. When considering the final set 
$A_{a-1}$, we also stop choosing paths for some $b \in B$ if the subcluster 
$V_b$ has size less than $\varepsilon m/4a$. Altogether this gives us a collection of 
close to $(k-1)(a-1)/c$ disjoint $\ell$-paths covering all but at most 
$\varepsilon m/2$ vertices in the clusters belonging to~$\F$. Doing this for all the copies 
of $\F_{k,\ell}$ in the $\F_{k,\ell}$-packing~$\mathcal{A}$ of $\R$ we obtain a 
collection of at most $|\mathcal{A}|(k-1)(a-1)/c\ll D$ disjoint $\ell$-paths covering all but at most 
$\varepsilon m/2$ vertices from each cluster, and hence all but at most 
$\varepsilon|\Hy|$ vertices of~$\Hy$, as required.
\endproof

\section{Proof of Theorem~\ref{main}}\label{sec:proof}

We shall use the following two results in our proof of Theorem~\ref{main}. The 
first says that if $1 \leq s \leq k-1$ and $\Hy$ is a large $k$-graph in which 
all sets of $s$ vertices have a large neighbourhood, then if we choose $R \subseteq 
V(\Hy)$ uniformly at random, with high probability all sets of $s$ vertices have 
a large neighbourhood in $R$.

\begin{lemma}[Reservoir Lemma] \label{reservoirlemma}
Suppose that $k \geq 2$, that $1 \leq s \leq k-1$, and that $1/n \ll \alpha, 
\mu, 1/k$. Let $\Hy$ be a $k$-graph of order $n$ with $d_\Hy(S) \geq \mu 
\binom{n}{k-s}$ for any set $S \in \binom{V(\Hy)}{s}$, and let~ $R$ be a subset 
of~$V(\Hy)$ of size $\alpha n$ chosen uniformly at random. Then the probability 
that $|N_\Hy(S) \cap \binom{R}{k-s}| \geq \mu \binom{\alpha n}{k-s} - 
n^{k-s-1/3}$ for every $S \in \binom{V(\Hy)}{s}$ is $1-o(1)$.
\end{lemma}

The proof of Lemma~\ref{reservoirlemma} is a standard probabilistic proof, which 
proceeds by applying Chernoff bounds to the size of the neighbourhood of each 
set $S$, and summing the probabilities of failure over all $S$. We omit the 
details.

The second result is the following theorem of Daykin and H\"aggkvist~\cite{DH},
giving an upper bound on the vertex degree needed to 
guarantee the existence of a perfect matching in a $k$-graph $\Hy$.

\begin{theo}[\cite{DH}] \label{perfmatching}
Suppose that $k \geq 2$ and $k|n$. Let $\Hy$ be 
a $k$-graph of order $n$ with minimum vertex degree at least 
$\frac{k-1}{k}\left(\binom{n-1}{k-1}-1\right)$. Then $\Hy$ contains a perfect matching.
\end{theo}

\removelastskip\penalty55\medskip\noindent{\bf Proof of Theorem~\ref{main}.}
In our proof we will use constants that satisfy the hierarchy 
$$1/n \ll 1/D \ll \eps \ll \alpha \ll c \ll \gamma \ll \gamma' \ll \eta \ll 
\eta' 
\ll 1/k.$$
Apply Lemma~\ref{absorb} to find an absorbing $\ell$-path $P_0$ in $\Hy$ which 
contains at most $\eta n/4$ vertices and which can absorb any set of at most 
$\alpha n$ $c$-good $(k-\ell)$-sets of vertices of $\Hy$. Define the 
$(k-\ell)$-graph $\G$ on the same vertex set as $\Hy$ to consist of all the 
$(k-\ell)$-sets of vertices of $\Hy$ which are $c$-good . Then by condition (1) 
of Lemma~\ref{absorb}, $d_\G(v) \geq \binom{n-1}{k-\ell-1} - \gamma n^{k-\ell-1} 
\geq (1-\gamma')\binom{n}{k-\ell-1}$ for every vertex $v$ in $V(\G)\sm V(P_0)$. 

Now, let $R$ be a set of $\alpha n$ vertices of $\Hy$ chosen uniformly at 
random. Then by Lemma~\ref{reservoirlemma}, with probability $1-o(1)$ we have 
that $|N_\G(v) \cap \binom{R}{k-\ell-1}| \geq (1-2\gamma') \binom{\alpha 
n}{k-\ell-1}$ for every vertex $v$ in $V(\G) \sm V(P_0)$. Likewise, with 
probability $1-o(1)$ we have that 
$$|N_\Hy(S) \cap R| \geq \left(\frac{1}{\lceil \frac{k}{k-\ell} \rceil 
(k-\ell)} + \frac{\eta}{2}\right)\alpha n.$$ 
for any $(k-1)$-set $S$ of vertices of $\Hy$. Finally, $\mathbb{E}[|R \cap 
V(P_0)|] = \alpha |P_0|$, and so with probability at least 1/2 we have that $|R 
\cap V(P_0)| \le \alpha \eta n/2$. Thus we may fix a choice of $R$ such that 
each of these three properties hold. Let $R' = R \sm V(P_0)$, so $|R'| \geq 
(1-\eta/2)\alpha n$. Then $|N_\G(v) \cap \binom{R'}{k-\ell-1}| \geq 
(1-\eta')\binom{\alpha n}{k-\ell-1}$ for every vertex $v$ in $V(\G) \sm V(P_0)$, 
and $|N_\Hy(S) \cap R'| \geq \frac{\alpha n}{\lceil \frac{k}{k-\ell} \rceil 
(k-\ell)}$ for any $(k-1)$-set $S$ of vertices of $\Hy$.

Let $V' = V(\Hy) \sm (V(P_0)\cup R)$, and let $\Hy' = \Hy[V']$ be the 
restriction of $\Hy$ to $V'$. Then as $|V(P_0)\cup R| \leq \eta n/2$, we must 
have
$$\delta(\Hy') \geq \left(\frac{1}{\lceil \frac{k}{k-\ell} \rceil 
(k-\ell)}+\frac{\eta}{2}\right) n.$$
We may therefore apply Lemma~\ref{pathcover} to $\Hy'$ to find a set of at most 
$D$ disjoint $\ell$-paths $P_1, \dots, P_q$ in $\Hy'$ which include all but at 
most $\eps n$ vertices of $\Hy'$. Let $X$ be the set of vertices not included in 
any of these $\ell$-paths, so $|X| \leq \eps n$.

For each $0 \leq i \leq q$, let $P_i^{beg}$ and $P_i^{end}$ be ordered ends of 
$P_i$. Next we shall find disjoint $\ell$-paths $P'_i$ for each $0 \leq i \leq 
q$, so that $P'_i$ is an $\ell$-path from
$P_i^{end}$ to $P_{i+1}^{beg}$ (where subindices are taken modulo $q+1$). The 
$\ell$-path $P'_i$ will only contain vertices from $R' \cup P_i^{end} \cup P_{i+1}^{beg}$, 
and will contain at most $8k^5$ vertices in total. So, suppose that we have 
found such $\ell$-paths $P'_0, \dots, P'_{i-1}$. Let
$R_i = (R' \cup P^{end}_i \cup 
P^{beg}_{i+1}) \sm  \bigcup_{j=0}^{i-1} V(P'_j)$. Then $\delta(\Hy[R_i]) \geq 
\frac{\alpha n}{\lceil \frac{k}{k-\ell}\rceil(k-\ell)} - 8k^5D \geq \alpha 
n/2k$, and so by Corollary~\ref{diameter} we can choose such an $\ell$-path 
$P'_i$ in~$\Hy[R_i]$. 

Then $C = P_0P_0'P_1P_1'\dots P_qP_q'$ is an $\ell$-cycle containing almost 
every vertex of $\Hy$. Indeed, $C$ contains every vertex of $\Hy$ except for 
those in $X$ and those in $R'$ not contained in any $P_i'$. So let $R'' = V(\Hy) 
\sm V(C)$. Then $(1-\eta)\alpha n\le |R''|\le(\alpha +\eps) n$. Since 
$(k-\ell)|n$ and $(k-\ell)\big||C|$ (as $C$ is an $\ell$-cycle), we also have 
$(k-\ell)\big||R''|$. Furthermore, $N_{\G[R'']}(v) \geq (1-2\eta')\binom{\alpha 
n}{k-\ell-1}$ for every vertex $v \in R''$. Since $k-\ell \geq 2$, 
Theorem~\ref{perfmatching} tells us that $\G[R'']$ contains a perfect matching, 
and so we can partition $R''$ into at 
most $\alpha n$ $c$-good $(k-\ell)$-sets of vertices of $\Hy$. Since $P_0$ can 
absorb any collection of at most $\alpha n$ $c$-good $(k-\ell)$-sets, there 
exists an $\ell$-path $Q_0$ with the same ordered ends as $P_0$ and such that 
$V(Q_0) = V(P_0) \cup R''$. 
Then $C' = Q_0P_0'P_1P_1'\dots P_qP_q'$ is a Hamilton $\ell$-cycle in $\Hy$, 
completing the proof of Theorem \ref{main}.

\medskip

\noindent
{\footnotesize
\noindent
Daniela K\"uhn, Richard Mycroft, Deryk Osthus, School of Mathematics,
University of Birmingham, Birmingham, B15 2TT, United Kingdom, 
{\{{\tt kuehn,mycroftr,osthus}\}\tt{@maths.bham.ac.uk}  }

\end{document}